\documentclass[12pt,english]{article} 
\usepackage[T2A]{fontenc}  
\usepackage{babel,amsfonts,latexsym,amssymb,amsmath}

\def \al{\alpha}
\def \be{\beta}
\def \ga{\gamma}
\def \dl{\delta}
\def \ep{\varepsilon}

\def \om{\omega}
\def \Ga{\Gamma}

\def \Si{\Sigma}

\def \Om{\Omega}

\def \operatorname#1{\mathop{\rm #1}}

\def\div{\operatorname{div}}
\def\osc{\operatorname{osc}}
\def\osc2{\operatorname{osc^2}}

\def\cd{\partial}

\def\rot{\operatorname{rot}}

\def\Q0{Q(x_0,t_0,R)}
\def\0{{x_0,t_0,R}}
\def\build#1_#2{\mathrel{\mathop{\kern 0pt#1}\limits_{#2}}}

\newtheorem{theorem}{Theorem}[section]

\newtheorem{corollary}{Corollary}[section]

\begin{document}

\title{On the Boundary Regularity of Weak Solutions to the MHD System}
\author{V.~Vyalov, T.~Shilkin\thanks{This work was supported by
RFBR-08-01-00372.}}

\date{December 2010}

\maketitle

\centerline{\it Dedicated to the jubilee of  Gregory A. Seregin}

\abstract{We prove the partial regularity of the boundary suitable weak solutions to the MHD system near the plane part of the boundary.}

\section{Introduction}

Assume $\Om\subset \mathbb R^3$ is a $C^2$-- smooth bounded domain and $Q_T=\Om\times (0,T)$. In this paper we investigate the boundary regularity of solutions to the principal system of magnetohydrodynamics (the MHD equations):
\begin{equation}
\left.
\begin{array}c
\partial_t v  + (v\cdot \nabla )v  - \Delta v +  \nabla
p = \rot H \times H
\\
 \div v =0
\end{array} \right\}
\quad\mbox{in } Q_T, \label{MHD_NS_000}
\end{equation}
\begin{equation}
\left.
\begin{array}c
\partial_t H    +\rot \rot H = \rot (v\times H)
\\ \div H =0
\end{array} \right\} \quad \mbox{in  }Q_T.
\label{MHD_Magnetic_000}
\end{equation}
Here  unknowns
are  the velocity field $v:Q_T\to \mathbb R^3$, pressure
$p:Q_T\to \mathbb R$, and the magnetic field $H:Q_T\to \mathbb
R^3$. We impose on $v$ and $H$ the initial and boundary conditions:
\begin{equation}
v|_{\cd\Om\times(0,T)}=0, \quad H_{\nu}|_{\cd\Om\times(0,T)}=0, \quad (\rot H)_{\tau}|_{\cd\Om\times(0,T)}=0,
\label{BC_000}
\end{equation}
\begin{equation}
v|_{t=0}=v_0, \qquad H|_{t=0}= H_0.
\label{IC_000}
\end{equation}
Here by $\nu$ we denote the outer normal to $\cd\Om$ and $H_{\nu}= H\cdot \nu$, $(\rot H)_\tau=\rot H- \nu(\rot H\cdot \nu)$.
We introduce the following definition:

\medskip\noindent{\bf Definition:} Assume $\Ga\subset \cd\Om$.
 The functions $(v,H, p)$ are called a {\it boundary suitable
 weak solution } to the system  (\ref{MHD_NS_000}), (\ref{MHD_Magnetic_000})
near $\Ga_T\equiv \Ga\times (0,T)$ if
\begin{itemize}
\item[1)] \ $v\in L_{2,\infty}(Q_T)\cap W^{1,0}_2(Q_T)\cap W^{2,1}_{\frac 98, \frac 32}(Q_T)$,

$H \in L_{2,\infty}(Q_T)\cap W^{1,0}_2(Q_T)$,
\item[2)] \ $p\in L_{\frac 32}(Q_T)\cap W^{1,0}_{\frac 98, \frac 32}(Q_T)$,
\item[3)] \ $\div v=0$, \ $\div H=0$ \ a.e. in \ $Q_T$,
\item[4)] \ $v|_{\cd \Om}=0$, \ $H_\nu|_{\cd\Om}=0$ \ in the sense of traces,
\item[5)] \ for any $w\in L_2(\Om)$ the functions
$$
t\mapsto \int\limits_{\Om} v(x,t)\cdot w(x)~dx \qquad\mbox{and}\qquad t\mapsto \int\limits_{\Om} H(x,t)\cdot w(x)~dx
$$
are continuous,

\item[6)] \ $(v,H)$ satisfy the following integral identities: for any $t\in [0,T]$
$$
\begin{array}c
\int\limits_{\Om} ~v(x,t)\cdot \eta(x,t)~dx -\int\limits_{\Om} ~v_0(x)\cdot \eta(x,0)~dx \ + \\ + \
\int\limits_0^t \int\limits_\Om ~ \Big( - v\cdot \cd_t \eta + (\nabla v  - v\otimes v + H \otimes H) : \nabla \eta - (p+\frac 12 |H|^2)\div \eta \Big)~dxdt \ = \ 0,
\end{array}
$$
for all \ $\eta \in W^{1,1}_{\frac 52}(Q_t)$ \ such that   \ $\eta|_{\cd\Om\times (0,t)}=0$,
$$
\begin{array}c
 \int\limits_{\Om}~ H(x,t)\cdot \psi (x,t)~dx  \ -  \  \int\limits_{\Om}~ H_0(x)\cdot \psi (x,0)~dx \ + \\ + \
\int\limits_0^t \int\limits_\Om ~ \Big( - H\cdot \cd_t \psi  + \rot H\cdot \rot \psi - (v\times H) \cdot \rot \psi \Big)~dxdt  \ = 0,
\end{array}
$$
for all \ $\psi \in W^{1,1}_{\frac 52}(Q_t)$ \ such that \  $\psi_\nu |_{\cd\Om\times (0,t)}=0$.

\item[7)] For every $z_0=(x_0,t_0)\in \Ga_T$  such that  $ \Om_R(x_0)\times (t_0-R^2, t_0)\subset Q_T$ and for any $\zeta\in C_0^\infty(B_R(x_0)\times(t_0-R^2,t_0])$ such that $\left.\frac{\cd\zeta}{\cd \nu}\right|_{\cd\Om}=0$ the following ``local energy inequality near $\Ga_T$'' holds:
\begin{equation}
\begin{array}c
\sup\limits_{t\in (t_0-R^2,t_0)}\int\limits_{\Om_R(x_0)}  \zeta \Big( |v|^2+|H|^2\Big)~dx \ + \\ + \ 2 \int\limits_{t_0-R^2}^{t_0}\int\limits_{\Om_R(x_0)} \zeta \Big( |\nabla v|^2+|\rot H|^2\Big)~dxdt \ \le \\
\le\  \int\limits_{t_0-R^2}^{t_0}\int\limits_{\Om_R(x_0)}  \Big(| v|^2 +|H|^2\Big) (\cd_t \zeta + \Delta \zeta )~dxdt \ + \\ + \
\int\limits_{t_0-R^2}^{t_0}\int\limits_{\Om_R(x_0)} \Big ( | v|^2+ 2\bar p \Big ) v\cdot\nabla \zeta~dxdt \ + \\
- \ 2 \int\limits_{t_0-R^2}^{t_0}\int\limits_{\Om_R(x_0)}   (H\otimes H):  \nabla^2 \zeta  ~dxdt \ + \\ + \
2 \int\limits_{t_0-R^2}^{t_0}\int\limits_{\Om_R(x_0)}  (v\times  H)(\nabla \zeta\times H)~dxdt
\end{array}
\label{LEI}
\end{equation}
We remark also that the following identity holds
$$
\begin{array}c
 \int\limits_{t_0-R^2}^{t_0}\int\limits_{\Om_R(x_0)}  (v\times  H)(\nabla \zeta\times H)~dxdt  \ =  \\ = \
\int\limits_{t_0-R^2}^{t_0}\int\limits_{\Om_R(x_0)} ~ (v \cdot \nabla \zeta) |H|^2  ~dxdt  \ -  \ \int\limits_{t_0-R^2}^{t_0}\int\limits_{\Om_R(x_0)} ~ (v \cdot H) (H\cdot \nabla \zeta)  ~dxdt
\end{array}
$$
\end{itemize}

\noindent
Here  $L_{s,l}(Q_T)$ is the anisotropic Lebesgue space equipped with the norm
$$
\|f\|_{L_{s,l}(Q_T)}:=
\Big(\int_0^T\Big(\int_\Om |f(x,t)|^s~dx\Big)^{l/s}dt\Big)^{1/l} ,
$$
and we use the following notation for the functional spaces:
$$
\gathered
W^{1,0}_{s,l}(Q_T)\equiv L_l(0,T; W^1_s(\Om))= \{ \ u\in L_{s,l}(Q_T): ~\nabla u \in L_{s,l}(Q_T) \ \},
\\
W^{2,1}_{s,l}(Q_T) = \{ \ u\in W^{1,0}_{s,l}(Q_T): ~\nabla^2 u, \ \cd_t u \in L_{s,l}(Q_T) \ \},
\\
\overset{\circ}{W}{^1_s}(\Om)=\{ \ u\in W^1_s(\Om):~ u|_{\cd\Om}=0 \ \},
\endgathered
$$
and the following notation for the norms:
$$
\gathered
\| u \|_{W^{1,0}_{s,l}(Q_T)}= \| u \|_{L_{s,l}(Q_T)}+ \|\nabla u\|_{L_{s,l}(Q_T)},  \\
\| u \|_{W^{2,1}_{s,l}(Q_T)}= \| u \|_{W^{1,0}_{s,l}(Q_T)}+ \| \nabla^2 u \|_{L_{s,l}(Q_T)}+\|\cd_t u\|_{L_{s,l}(Q_T)},  \\
\endgathered
$$

\begin{theorem} \label{Existense_Suitable}
For any sufficiently smooth divergent-free \ $v_0$,  $H_0$ satisfying (\ref{BC_000})  there exists at least one boundary suitable weak solution near $\cd\Om\times (0,T)$ which satisfies the initial conditions :
$$\| v(\cdot, t) -v_0(\cdot)\|_{L_2(\Om)}\to 0, \quad \| H(\cdot, t) -H_0(\cdot)\|_{L_2(\Om)}\to 0 \quad\mbox{as} \quad t\to+0,$$
and additionally satisfies the global energy inequality
$$
\gathered
\| v\|_{L_{2,\infty}(Q_T)}+ \| H\|_{L_{2,\infty}(Q_T)} + \| \nabla v\|_{L_{2}(Q_T)} + \| \rot H\|_{L_{2}(Q_T)} \le \\ \le \| v_0\|_{L_2(\Om)}+ \| H_0\|_{L_2(\Om)}
\endgathered
$$
\end{theorem}

The global existence of weak solutions to the MHD equations (\ref{MHD_NS_000}) --- (\ref{BC_000}) was established originally  in \cite{LadSol}.
We  sketch the proof of Theorem \ref{Existense_Suitable} in the Section 3.

\section{Main Results}

\bigskip
Denote $B(x_0,r)$ the open ball in $\mathbb R^3$ of radius $r$ centered at $x_0$ and denote by $B^+(x_0,r)$ the half--ball $\{ x\in B(x_0,r)~|~ x_3>0\}$. For $z_0=(x_0,t_0)$ denote $Q(z_0,r)=B(x_0,r)\times (t_0-r^2, t_0)$, $Q^+(z_0,r)=B^+(x_0,r)\times (t_0-r^2, t_0)$. In this paper we shall use the abbreviations: $B(r)= B(0,r)$, $B^+(r)= B^+(0,r)$ etc, $B=B(1)$, $B^+=B^+(1)$ etc.

\medskip
The main results of the present paper are the following theorems:

\begin{theorem}\label{Fixed_r}
There exists an absolute constant $\ep_*>0$ with the following property. Assume $(v,H,p)$ is a boundary suitable weak solution in $Q_T$ and assume $z_0=(x_0,t_0)\in \cd\Om\times (0,T)$ is such that $x_0$ belongs to the plane part of $\cd\Om$. If there exists $r_0>0$ such that $Q^+(z_0,r_0)\subset Q_T$ and
$$
\frac 1{r_0^2} \int\limits_{Q^+(z_0,r_0)} \Big(~ |v|^3+|H|^3+|p|^{\frac 32}~\Big) dxdt \  < \ \ep_*,
$$
then the functions $v$ and $H$ are H\" older continuous on the closure of $ Q^+(z_0, \frac {r_0}2)$.
\end{theorem}


\begin{theorem}\label{CKN_theorem}
For any  $K>0$ there exists $\ep_0(K)>0$ with the following property.
Assume $(v,H,p)$ is a boundary suitable weak solution in $Q_T$ and assume $z_0=(x_0,t_0)\in \cd\Om\times (0,T)$ is such that $x_0$ belongs to the plane part of $\cd\Om$.
If
\begin{equation}
\limsup\limits_{r\to 0}\Big(~\frac 1r \int\limits_{Q(z_0,r)} |\nabla H|^2~dxdt ~\Big)^{1/2}\ < \ K
\label{ep-regularity-1}
\end{equation}
and
\begin{equation}
\limsup\limits_{r\to 0} \Big(~\frac 1r \int\limits_{Q(z_0,r)} |\nabla v|^2~dxdt~\Big)^{1/2} \ < \ \ep_0,
\label{ep-regularity-2}
\end{equation}
then there exists $\rho_*>0$ such that the functions $v$ and $H$ are H\" older continuous on the closure of $ Q^+(z_0,\rho_*)$.

\end{theorem}

\begin{theorem}\label{Partial_Regularity}
Assume $(v,H,p)$ is a boundary suitable weak solution in $Q_T$ and denote by $\Ga$ the plane part of $\cd \Om$. Then there exists a closed set $\Si\subset \Ga$ such that
for any $z_0\in (\Ga\setminus \Si)\times (0,T]$ the functions $(v,H)$ are H\" older continuous in some neighborhood of $z_0$, and
$$
\mathcal P^1(\Si) \ = \ 0,
$$
where $\mathcal P^1(\Si)$ is the one-dimensional parabolic Hausdorff measure of $\Si$.
\end{theorem}

Our Theorem \ref{Partial_Regularity} presents for the MHD equations a result which is a boundary analogue of  the famous Caffarelli--Kohn--Nirenberg (CKN)
 theorem  for the Navier-Stokes system, see \cite{CKN}, see also \cite{Lin}.
 The boundary regularity of solutions to the Navier-Stokes equations was originally investigated  by G.~Seregin in \cite{Seregin_JMFM} and \cite{Seregin_Aa}  in the case of a plane part of the boundary and by G.~Seregin, T.~Shilkin, and V.~Solonnikov in \cite{SSS} in the case of a curved boundary. The internal partial regularity of  solutions to the MHD system was originally proved by C.~He and Z.~Xin in \cite{China}, see also \cite{Vya}, \cite{Vya1}. Note that thought using the methods of our paper one can prove various $\ep$--regularity conditions involving various scale--invariant  functionals (such as in \cite{China}), in the present paper    we concentrate on the condition (\ref{ep-regularity-1}), (\ref{ep-regularity-2}) as this condition provides the optimal estimate of the Hausdorff measure of the singular set $\Si$ in Theorem \ref{Partial_Regularity}.

\medskip
Our paper is organized as follows: in Section 3 we outline the proof of Theorem \ref{Existense_Suitable}. In Section 4 we prove that the weak solutions of the linearized MHD equations are H\" older continuous up to the boundary.  Section 5 contains the proof of the Decay Lemma and the sketch of the proof of Theorem \ref{Fixed_r}.  Section 6 is concerned with the estimate of some Morrey functional for weak solutions to the heat equation near the boundary. These estimates together with the estimates of the scale invariant  energy functionals obtained in Section 7 turn to be crucial for the key estimate (\ref{Estimate_varE}) in the Section 8. Finally, in Section 8 we present the proofs of Theorems \ref{CKN_theorem} and \ref{Partial_Regularity}.

\vskip15pt
\section{Local  Energy Inequality  and  \\ Existence of Boundary Suitable \\ Weak Solutions}
\setcounter{equation}{0}

\medskip
\bigskip\noindent
Denote by $\om$ the usual smooth Sobolev kernel and  $\om_\dl(x)=\dl^{-3}\om(\dl x)$. Denote by $\om_\dl* H $ the convolution of $H$ with $\om_\dl$.
Assume  $v_0^\dl$ and  $H^\dl_0$ are sooth divergent--free functions satisfying the boundary conditions (\ref{BC_000}) and additionally satisfying  the properties
$$
v_0^\dl\to v_0, \quad H_0^\dl\to H_0 \quad \mbox{in}\quad L_2(\Om) \qquad\mbox{as}\quad\dl\to0.
$$
Consider the problem
\begin{equation}
\left.
\gathered
\partial_t v +  ((\om_\dl * v)\cdot \nabla )v - \Delta v + \nabla p = \rot H \times (\om_{\dl}* H) \\
\div v = 0
\endgathered
\quad
\right\} \qquad\mbox{in} \quad Q_T,
\label{deltaeq1}
\end{equation}
\begin{equation}
\left.
\gathered
 \partial_t H + \rot\rot H = \rot ( v \times (\om_{\dl}*H)) \\
 \div H =0
 \endgathered\quad
\right\} \qquad\mbox{in} \quad Q_T,
\label{deltaeq2}
\end{equation}
$$
\gathered
v|_{\cd\Om\times (0,T)}=0, \qquad H_{\nu}|_{\cd\Om\times (0,T)}=0, \qquad (\rot H)_\tau |_{\cd\Om\times (0,T)}=0, \\
v|_{t=0}=v_0^\dl, \qquad H|_{t=0}=H_0^\dl.
\endgathered
$$
Using Galerkin approximations and energy estimates it is not difficult to prove existence of the unique smooth solution $(v^\dl, H^\dl, p^\dl)$ to this problem.
Moreover, these functions satisfy the global energy inequality
$$
\gathered
\| v^\dl\|_{L_{2,\infty}(Q_T)}+ \| H^\dl\|_{L_{2,\infty}(Q_T)} + \| v^\dl\|_{W^{1,0}_{2}(Q_T)} + \|  H^\dl\|_{W^{1,0}_{2}(Q_T)} \le \\ \le c~\Big(~\| v_0\|_{L_2(\Om)}+ \| H_0\|_{L_2(\Om)}~\Big)
\endgathered
$$
From this estimate using the interpolation inequality  we obtain the estimate $$\| ((\om_\dl * v^\dl)\cdot \nabla )v^\dl\|_{L_{\frac 98, \frac 32}(Q_T)}+\| \rot H^\dl \times (\om_{\dl}* H^\dl)\|_{L_{\frac 98, \frac 32}(Q_T)} \ \le \ c.$$
Applying the coercive estimates for the linear Stokes problem (see, for example, \cite{Solonnikov_Uspekhi}), we obtain
$$
\| v^\dl\|_{W^{2,1}_{\frac 98, \frac 32}(Q_T)} + \| \nabla p^\dl \|_{L_{\frac 98, \frac 32}(Q_T)} \ \le \ c.
$$
So, the only thing we need is to verify the local energy inequality (\ref{LEI}).  Then the result of Theorem \ref{Existense_Suitable} follows if we pass to the limit as $\dl\to 0$.

To simplify notations below we omit the index $\dl$ and denote by $(v, H, p)$ the smooth functions $(v^\dl,H^\dl, p^\dl)$. Take $z_0=(x_0, t_0)\in \cd\Om\times (0,T)$ and choose $R$ so that $\Om_R(x_0)\times (t_0-R^2, t_0)\subset Q_T$. Without loss of generality we can put $x_0=0$, $t_0=0$.
Assume $\zeta\in C_0^\infty (B_R\times (-R^2,0])$ satisfies $\frac{\cd \zeta}{\cd\nu}\big|_{\cd\Om}=0$ and multiply the equation (\ref{deltaeq1})
by the test--function  $\eta=\zeta v$. Integrating the result over $\Om$, integrating by parts and taking into account the relation
$$
\begin{array}c
\int\limits_{\Om} ~ (\om_\dl*v_i) v_{j,i} \zeta v_j ~dx \ = \ - \frac 12~ \int\limits_{\Om} ~ |v|^2 (\om_\dl*v) \cdot \nabla \zeta  ~dx .
\end{array}
$$
 we obtain
\begin{equation}
\begin{array}c
\frac 12 \frac {d}{dt}~ \int\limits_{\Om}  ~\zeta |v|^2 \ + \
\int\limits_{\Om}
\zeta |\nabla v |^2 ~dx \ = \   \frac 12 \int\limits_{\Om}  |v|^2 (\cd_t\zeta + \Delta \zeta )~dx  \ +  \\ + \  \int\limits_{\Om} ~ \bar p~ v \cdot \nabla \zeta  ~dx \ + \ \frac 12~ \int\limits_{\Om} ~ |v|^2 (\om_\dl*v) \cdot \nabla \zeta  ~dx \ + \\ +
  \  \int\limits_{\Om} ~ (\rot H\times (\om_\dl*H))\cdot \zeta v~dx = 0 .
\end{array}
\label{Lei for v}
\end{equation}
Now multiply the equation for the magnetic field by $\zeta H$. Integrating by parts and taking into account the boundary conditions (\ref{BC_000}) after routine calculations we arrive at
$$
\begin{array}c
\frac 12 \frac {d}{dt}~ \int\limits_{\Om}  ~\zeta |H|^2~dx \ + \
\int\limits_{\Om}
\zeta |\rot H |^2 ~dx \ = \\ = \   \frac 12 \int\limits_{\Om}  |H|^2 (\cd_t\zeta + \Delta \zeta )~dx  \ - \  \int\limits_{\Om} ~ H\otimes H : \nabla^2\zeta  ~dx \ + \\ +
  \  \int\limits_{\Om} ~ \rot (v\times (\om_\dl * H))\cdot \zeta H ~dx.
\end{array}
$$
Integrating by parts, we obtain
$$
\begin{array}c
\int\limits_{\Om} ~ \rot (v\times (\om_\dl* H))\cdot \zeta H ~dx \ = \  \int\limits_{\Om} ~ \zeta (v\times (\om_\dl *H))\cdot \rot  H ~dx  \ + \\ + \ \int\limits_{\Om} ~ (v\times (\om_\dl* H))\cdot (\nabla \zeta \times H) ~dx
\end{array}
$$
Note that
$$
\begin{array}c
\int\limits_{\Om} ~ \zeta (v\times (\om_\dl *H))\cdot \rot  H ~dx \ =  \ \int\limits_{\Om} ~ \zeta  ( (\om_\dl *H) \times \rot H) \cdot v ~dx
\end{array}
$$
Hence we obtain
\begin{equation}
\begin{array}c
\frac 12 \frac {d}{dt}~ \int\limits_{\Om}  ~\zeta |H|^2~dx \ + \
\int\limits_{\Om}
\zeta |\rot H |^2 ~dx \ =  \   \frac 12 \int\limits_{\Om}  |H|^2 (\cd_t\zeta + \Delta \zeta )~dx  \ - \\ - \  \int\limits_{\Om} ~ H\otimes H : \nabla^2\zeta  ~dx \ +
  \  \int\limits_{\Om} ~ ((\om_\dl *H)\times \rot H)\cdot \zeta v ~dx \ + \\ + \  \int\limits_{\Om} ~ (v\times (\om_\dl *H))\cdot (\nabla \zeta \times H) ~dx.
\end{array}
\label{Lei for H}
\end{equation}
Adding  identities (\ref{Lei for v}) and (\ref{Lei for H}) together,  using the  conciliation
$$
\begin{array}c
\int\limits_{\Om} ~ (\rot H \times (\om_\dl * H))\cdot \zeta H ~dx \ + \ \int\limits_{\Om} ~ ((\om_\dl *H)\times \rot H)\cdot \zeta v ~dx = 0
\end{array}
$$
and passing to the limit as $\dl\to 0$ we arrive at (\ref{LEI}).  Theorem \ref{Existense_Suitable} is proved.

\bigskip
In conclusion of this section we introduce the following version of the local energy inequality near the plane part of the boundary.

\begin{theorem}\label{LEI near plane}
Assume  $(v,H,p)$ is a boundary suitable weak solution satisfied the MHD equations in a domain containing the set $Q^+=B^+\times (-1,0)$ such that the plane part of $Q^+$ belongs to the boundary $\cd\Om\times (-1,0)$. Obviously then we have
$$
v|_{x_3=0}=0, \quad H_3|_{x_3=0}=0, \quad H_{1,3}|_{x_3=0}= H_{2,3}|_{x_3=0}=0.
$$
Let $\zeta\in C_0^\infty(B\times (-1,0])$ be a cut-off function such that $\zeta_{, 3}|_{x_3=0}=0$.
Assume $b\in \mathbb R^3$ is an arbitrary constant vector of the form $b=(b_1, b_2, 0)$. Then the following inequality holds
\begin{equation}
\begin{array}c
\sup\limits_{t\in (-1,0)}\int\limits_{B^+}  \zeta \Big( |v|^2+|\bar H|^2\Big)~dx \ +  \ 2 \int\limits_{Q^+} \zeta \Big( |\nabla v|^2+|\rot\bar H|^2\Big)~dz \ \le \\
\le\  \int\limits_{Q^+}  \Big(| v|^2 +|\bar H|^2\Big) (\cd_t \zeta + \Delta \zeta )~dz \ + \
\int\limits_{Q^+} \Big ( | v|^2+ 2\bar p \Big ) v\cdot\nabla \zeta~dz \ + \\
- \ 2 \int\limits_{Q^+}   (\bar H\otimes \bar H):  \nabla^2 \zeta  ~dz \ + \
2 \int\limits_{Q^+}  (v\times   H)(\nabla \zeta\times \bar H)~dz
\end{array}
\label{LEI_near_plane_part}
\end{equation}
where $\bar H= H-b$.
\end{theorem}

\noindent
{\bf Proof:} The relation (\ref{LEI_near_plane_part}) is a combination of (\ref{LEI}) with he relation obtained from (\ref{MHD_Magnetic_000}) multiplied by the  test function $\psi=\zeta b$ where $b=(b_1, b_2, 0) $ is a constant vector. The proof is routine and we omit it.

\section{Linear Estimate}
\setcounter{equation}{0}

\bigskip
\begin{theorem}
\label{Linear_Estimate_Theorem}
For any $M>0$ there exists $c(M)>0$ such that for any $a \in \mathbb R^3$ such that $a=(a_1, a_2, 0)^T$ and  $|a|\le M$, and for any $(u, h, q)$ satisfying the linear system
\begin{equation}
\gathered
\partial_t u   - \Delta u +  \nabla q = \rot h \times   a
\\
 \div u =0
\endgathered  \label{MHD_NS_0}
\end{equation}
\begin{equation}
\gathered
\partial_t h     -  \Delta h = \rot (u\times  a)
\\ \div h =0
\endgathered
\label{MHD_Magnetic_0}
\end{equation}
\begin{equation}
\gathered
u|_{x_3=0}=0, \\
\begin{array}c
h_3|_{x_3=0}=0, \qquad \frac{\cd h_\al}{\cd x_3}\big|_{x_3=0}\  = \ 0, \qquad \al=1,2,
\end{array}
\endgathered
\label{BC_0}
\end{equation}
the following estimate holds:
\begin{equation}
\gathered
\| u \|_{C^{\frac 13, \frac 16}(Q^+(\frac 12))} + \| h \|_{C^{\frac 13, \frac 16 }(Q^+(\frac 12))}   \ \le  \\
\le \ c(M)\Big( \| u \|_{L_3(Q^+)} + \| h -b\|_{L_3(Q^+)} + \| q -c\|_{L_{3/2}(Q^+)}  \Big)
\endgathered
\label{Linear_Estimate}
\end{equation}
Here  $b\in \mathbb R^3$ is an arbitrary  vector of the form $b = (b_1, b_2, 0)^T$, and $c\in \mathbb R$ is an arbitrary constant.

\end{theorem}

\noindent
{\bf Proof:}

\medskip
\noindent
{\bf 1.} Similar to Theorem \ref{LEI near plane} we obtain
\begin{equation}
\gathered
\| u \|_{L_{2,\infty}(Q^+(\frac 9{10}))} \ + \ \| h \|_{L_{2,\infty}(Q^+(\frac 9{10}))}  \ + \\ + \
\| \nabla u \|_{L_{2}(Q^+(\frac 9{10}))} \ + \ \| \nabla h \|_{L_{2}(Q^+(\frac 9{10}))} \ \le \\
\le \ c(M) ~\Big( \| u\|_{L_2(Q^+)} + \| h-b\|_{L_2(Q^+)} \ + \ \| u \bar q  \|_{L_1(Q^+)}\Big)
\endgathered
\label{Energy}
\end{equation}
Here $\bar q\equiv q-c$. Note that the right-hand side of the last inequality can be estimated by the right-hand side of (\ref{Linear_Estimate}) via H\" older inequality.

\medskip
\noindent
{\bf 2.} For the function $u$ satisfying the system (\ref{MHD_Magnetic_0}) we have the following estimate with arbitrary $s$, $l\in (1,+\infty)$ and $0< \rho <r <\frac 9{10}$ (see \cite{Seregin_ZNS271}):
\begin{equation}
\gathered
\| u \|_{W^{2,1}_{s, l}(Q^+(\rho))} \ \le    \ c ~\Big(\| \rot h\|_{L_{s, l}(Q^+(r))} \ +  \| u \|_{L_{3}(Q^+)} + \| \bar q  \|_{L_{3/2}(Q^+)}\Big)
\endgathered
\label{Anisotropic_Estimate}
\end{equation}
Here $C$ depends on $M$, $r$, $\rho$, $s$, and $l$.

\medskip
\noindent
{\bf 3.} Denote by $h^*$ the following extension of the function $h$ from $Q^+$ onto $Q$:  we take $h^*(x,t)=h(x,t)$ for $x_3\ge 0$ and take
$$
\left\{ \
\gathered
h_\al^*(x_1, x_2, x_3,t) = \ ~ h_\al(x_1, x_2, - x_3,t) \\
\ h_3^*(x_1, x_2, x_3,t) = -h_3(x_1, x_2, - x_3,t)
\endgathered\right.
\qquad \mbox{for } \ x_3<0.
$$
Denote  $ g =\rot  (u\times H)$ and let $g^*$ be the following extension on $g$ from $Q^+$ onto $Q$: the components $g_{\al}^*$, $\al=1,2$ are the even extensions of $g_{\al}$, and the component $g_3^*$ is the odd extension of $g_3$. Then the functions $h^*$, $g^*$ satisfy the equation
\begin{equation}
\gathered
\partial_t h^*     -  \Delta h^* = g^*\qquad \mbox{in } \ Q.
\endgathered
\label{Heat_equation}
\end{equation}

\medskip
\noindent
{\bf 4.} For the function $h^*$ satisfying the heat equation (\ref{Heat_equation}) the  estimate similar to (\ref{Anisotropic_Estimate}) holds.
$$
\gathered
\| h^* \|_{W^{2,1}_{s, l}(Q(\rho))} \ \le    \ C ~\Big(\| g^* \|_{L_{s, l}(Q(r))} \ +  \| h^*-b \|_{L_{3}(Q)} \Big)
\endgathered
$$
Here $s$, $l\in (1,+\infty)$ and $0< \rho <r \le\frac 9{10}$ are arbitrary and  $C$ depends on $M$, $r$, $\rho$, $s$, and $l$. This estimate provides the inequality
\begin{equation}
\gathered
\| h \|_{W^{2,1}_{s, l}(Q^+(\rho))} \ \le    \ C ~\Big(\| \nabla u \|_{L_{s, l}(Q^+(r))} \ +  \| h -b \|_{L_{3}(Q^+)} \Big)
\endgathered
\label{Anisotropic_Estimate_2}
\end{equation}

\medskip
\noindent
{\bf 5.} First we apply (\ref{Anisotropic_Estimate_2}) with $s=l=2$ and $\rho=\frac 45$, $r =\frac 9{10} $. Taking into account the energy estimate (\ref{Energy}) we obtain the estimate of the norm $\| h\|_{W^{2,1}_2(Q^+(\frac 45))}$ by the right-hand side of (\ref{Linear_Estimate}). Then using the imbedding theorem $W^{2,1}_{2}(Q^+(\frac 45))\hookrightarrow W^{1,0}_{\frac {10}3}(Q^+(\frac 45))$ we obtain the estimate of the norm $\| \nabla h\|_{L_{\frac {10}3} (Q^+(\frac 45))}$.

\medskip
\noindent
{\bf 6.} Now applying (\ref{Anisotropic_Estimate}) with $s=\frac{10}3$, $l=\frac 32$ and $\rho=\frac 7{10}$, $r =\frac 4{5} $ we obtain the estimate
$$
\| u \|_{W^{2,1}_{\frac {10}3, \frac 32}(Q^+(\frac {7}{10}))} \ \le    \ C ~\Big(\| \rot h\|_{L_{\frac {10}3, \frac 32}(Q^+(\frac 45))} \ +  \| u \|_{L_{3}(Q^+)} + \| \bar q  \|_{L_{3/2}(Q^+)}\Big)
$$
By H\" older inequality we estimate $\| \rot h\|_{L_{\frac {10}3, \frac 32}(Q^+(\frac 45))}$ by the norm $\| \nabla h\|_{L_{\frac {10}3} (Q^+(\frac 45))}$ which was already estimated on the previous step. On the other hand, by the imbedding theorem we obtain
\begin{equation}
\| \nabla u \|_{L_{\frac 32}(-\left(\frac {7}{10}\right)^2, 0 ; L_\infty (B^+(\frac 7{10})))} \le C\| u \|_{W^{2,1}_{\frac {10}3, \frac 32}(Q^+(\frac {7}{10}))}
\label{nabla_u_estimate}
\end{equation}

\medskip
\noindent
{\bf 7.}  The estimate (\ref{nabla_u_estimate}) implies in particular that $\nabla u \in L_{9,\frac 32}(Q^+(\frac 7{10}))$.
Applying (\ref{Anisotropic_Estimate_2}) with $s=2$, $l=\frac 32$ we obtain the estimate
$$
\| h\|_{W^{2,1}_{9, \frac 32}(Q^+(\frac 35))} \ \le \ C~ \Big( \| \nabla u \|_{L_{9,\frac 32}(Q^+(\frac 7{10}))} + \| h-b\|_{L_3(Q^+)}\Big)
$$

\medskip
\noindent
{\bf 8.}  Finally, we apply (\ref{Anisotropic_Estimate}) with $s=9$, $l=\frac 32$ and obtain
$$
\gathered
\| u \|_{W^{2,1}_{9, \frac 32}(Q^+(\frac 12))} \ \le    \ C ~\Big(\| \rot h\|_{L_{9, \frac 32}(Q^+(\frac 35))} \ +  \| u \|_{L_{3}(Q^+)} + \| \bar q  \|_{L_{3/2}(Q^+)}\Big)
\endgathered
$$

\medskip
\noindent
{\bf 9.} Gathering all the estimates together we arrive at the estimate
$$
\gathered
\| u \|_{W^{2,1}_{9, \frac 32}(Q^+(\frac 12))}  + \| h \|_{W^{2,1}_{9, \frac 32}(Q^+(\frac 12))}  \ \le \\ \le \
C(M)\Big( \| u \|_{L_3(Q^+)} + \| h -b\|_{L_3(Q^+)} + \| q -c\|_{L_{3/2}(Q^+)}  \Big)
\endgathered
$$
The statement of the theorem follows now from the  imbedding $W^{2,1}_{s,l}(Q_T)\hookrightarrow C^{\al,\frac \al2}(\bar Q_T)$ as $\al=2-\frac 3s -\frac 2l>0$,
where $s=9$, $l=\frac 32$, and $\al=\frac 13$. Theorem \ref{Linear_Estimate_Theorem} is proved.

\medskip
\begin{corollary}\label{Corollary_Linear_Estimate}
Let us introduce the functional
\begin{equation}
\gathered
Y_\tau (v,H,p) \ : =  \ \Big(\ - \!\!\!\!\!\!\!\!\!\!\!\!  \ \int\limits_{Q^+(\tau)} ~|v|^3~dxdt\ \Big)^{1/3} \ + \\ \ \qquad\qquad + \
\tau ~\Big(\ - \!\!\!\!\!\!\!\!\!\!\!\!  \ \int\limits_{Q^+(\tau)} ~|p-[p]_{B^+(\tau)}|^{3/2}~dxdt\ \Big)^{2/3} \ + \\
 \ \qquad\qquad + \ \Big(\ - \!\!\!\!\!\!\!\!\!\!\!\!  \ \int\limits_{Q^+(\tau)} ~|H-b_\tau (H)|^3~dxdt\ \Big)^{1/3}
\endgathered
\label{Def_Y}
\end{equation}
where \quad $b_\tau(H):= ( (H_1)_{Q^+(\tau)}, (H_2)_{Q^+(\tau)}, 0)$. Then for any $M>0$ there is a constant $C_l(M)>0$ such that for any solution $(u,h,q)$ of the linear system (\ref{MHD_NS_0}), (\ref{MHD_Magnetic_0}), (\ref{BC_0}), the following estimate holds:
$$
Y_\tau (u,h,q) \ \le C_l(M) ~\tau^{1/3} ~Y_1(u,h,q).
$$

\end{corollary}

\newpage

\section{The Decay Estimate and \\ the Proof of Theorem \ref{Fixed_r}}
\setcounter{equation}{0}

\bigskip\medskip
In this section we consider the MHD system
\begin{equation}
\left.
\gathered
\partial_t v + (v\cdot \nabla )v  - \Delta v +  \nabla
p = \rot H \times H
\\
 \div v =0
\endgathered  \quad \right\}
\quad\mbox{in } Q^+, \label{MHD_NS}
\end{equation}
\begin{equation}
\left.
\gathered
\partial_t H    -  \Delta H = \rot (v\times H)
\\ \div H =0
\endgathered \quad \right\} \quad \mbox{in  }Q^+,
\label{MHD_Magnetic}
\end{equation}

\begin{equation}
\left.
\gathered
v|_{x_3=0}=0, \\
\begin{array}c
H_3|_{x_3=0}=0, \qquad \frac{\cd H_1}{\cd x_3}\big|_{x_3=0}\ =\ \frac{\cd H_2}{\cd x_3}\big|_{x_3=0} \ = \ 0
\end{array}
\endgathered \quad \right\}
\label{BC}
\end{equation}

\medskip\noindent
We denote by $Y_\tau (v,H,p)$ the functional introduced in (\ref{Def_Y}). We also denote by $Y_\tau (v)$, $\tilde Y_\tau (H)$, and $\hat Y_\tau (p)$ the functionals
$$
\gathered
Y_\tau (v) \ : =  \ \Big(\ - \!\!\!\!\!\!\!\!\!\!\!\!  \ \int\limits_{Q^+(\tau)} ~|v|^3~dxdt\ \Big)^{1/3}, \\
\tilde Y_\tau (H) \ := \ \Big(\ - \!\!\!\!\!\!\!\!\!\!\!  \ \int\limits_{Q^+(\tau)} ~|H-b_\tau (H)|^3~dxdt\ \Big)^{1/3}, \\
\hat Y_\tau (p)\ := \
\tau ~\Big(\ - \!\!\!\!\!\!\!\!\!\!\!  \ \int\limits_{Q^+(\tau)} ~|p-[p]_{B^+(\tau)}|^{3/2}~dxdt\ \Big)^{2/3}
\endgathered
$$

\begin{theorem} \label{Decay estimate theorem}
There exists an absolute constant $\ep_0>0$ such that for any $M>0$ there exists $C_*=C_*(M)$ with the following properties.  For any boundary suitable weak solution $(v,H,p)$ of the MHD system (\ref{MHD_NS}), (\ref{MHD_Magnetic}), (\ref{BC}), the following implication holds:
if
$$
Y_1(v,H,p) \ < \ \ep_0
$$
and
$$
|(H_1)_{Q^+}|+|(H_2)_{Q^+}|\ \le \  M
$$
then
\begin{equation}
Y_\tau (v,H,p) \ \le \ C_* ~\tau^{1/3}~Y_1 (v,H,p)
\label{Decay_Estimate}
\end{equation}
\end{theorem}

\medskip
\noindent
{\bf Proof:}

\medskip
\noindent
{\bf 1.} Arguing by contradiction we assume there exists a sequence of numbers $\ep_m\to 0$, and a sequence of boundary suitable weak solutions $(v^m,H^m, p^m)$ such that
$$
Y_1(v^m,H^m,p^m) \ = \ \ep_m \ \to \ 0
$$
and
$$
Y_\tau (v^m, H^m, p^m)\ \ge \  C_* \tau^{1/3} \ep_m
$$

\medskip
\noindent
{\bf 2.}
Let us introduce functions
$$
\gathered
u^m(x, t) \ = \  \frac{1}{\ep_m}~ v^m(x, t), \\
q^m(x, t) \ = \ \frac{1}{\ep_m}~ \Big(p^m(x, t) - [p^m]_{B^+}(t)\Big), \\
h^m(x,t) \ = \  \frac{1}{\ep_m}~ \Big(H^m(x, t) - a^m\Big),\\
a^m \ = \ b_1(H^m)
\endgathered
$$
Then
\begin{equation}
Y_1(u^m,h^m,q^m)\ = \  1,
\qquad
Y_\tau (u^m, h^m, q^m )\ge C_* \tau^{1/3}
\label{Contradiction_asumptions}
\end{equation}
and $(u^m, h^m, q^m)$ satisfy the following equations in $\mathcal D'(Q^+)$
\begin{equation}
\gathered
\partial_t u^m + \ep_m (u^m\cdot \nabla )u^m  - \Delta u^m +  \nabla
q^m = \rot h^m \times ( \ep_m h^m +  a^m)
\\
 \div u^m =0
\endgathered  \label{MHD_NS_m}
\end{equation}
\begin{equation}
\gathered
\partial_t h^m     -  \Delta h^m = \rot \big(u^m\times (\ep_m h^m+ a^m)\big)
\\ \div h^m =0
\endgathered
\label{MHD_Magnetic_m}
\end{equation}
\begin{equation}
\left.
\gathered
u^m|_{x_3=0}=0, \\
\begin{array}c
h^m_3|_{x_3=0}=0, \qquad \frac{\cd h^m_1}{\cd x_3}\big|_{x_3=0}\ =\ \frac{\cd h^m_2}{\cd x_3}\big|_{x_3=0} \ = \ 0
\end{array}
\endgathered \quad \right\}
\label{BC_m}
\end{equation}

\medskip
\noindent
{\bf 3.} The conditions (\ref{Contradiction_asumptions}) imply in particular the boundedness
\begin{equation}
\sup\limits_m~\Big(\| u^m\|_{L_3(Q^+)} + \| h^m \|_{L_3(Q^+)} + \| q^m \|_{L_{\frac 32}(Q^+)}\Big) \ < \ +\infty
\label{Basic_boundedness}
\end{equation}
From the local energy inequality near the boundary (\ref{LEI_near_plane_part}) we obtain
\begin{equation}
\gathered
\| u^m \|_{L_{2,\infty}(Q^+(\frac 9{10}))} + \| h^m \|_{L_{2,\infty}(Q^+(\frac 9{10}))} \ + \\  + \  \| u^m \|_{W^{1,0}_{2}(Q^+(\frac 9{10}))} + \| h^m \|_{W^{1,0}_{2}(Q^+(\frac 9{10}))} \ \le C.
\endgathered
\label{Gradient boundedness}
\end{equation}
From the equations (\ref{MHD_NS_m}), (\ref{MHD_Magnetic_m}), (\ref{BC_m}) we also obtain the estimate
$$
\| \cd_t u^m \|_{L_{\frac 53}(-1, 0; W^{-1}_{\frac 53} (B^+)) } + \| \cd_t h^m \|_{L_{\frac 53}(-1, 0; W^{-1}_{\frac 53} (B^+)) }\ \le \ C.
$$

\medskip
\noindent
{\bf 4.}
Hence we can extract subsequences
\begin{equation}
\begin{array}c
u^m \ \rightharpoonup \ u \quad \mbox{in}  \quad L_3(Q^+), \\
h^m \ \rightharpoonup \ h \quad \mbox{in}  \quad L_3(Q^+), \\
q^m \ \rightharpoonup \ q \quad \mbox{in}  \quad L_{\frac 32}(Q^+), \\
\end{array}
\label{Weak_Convergence0}
\end{equation}
\begin{equation}
\begin{array}c
u^m \ \rightharpoonup \ u \quad \mbox{in}  \quad W^{1,0}_2(Q^+(\frac 9{10})), \\
h^m \ \rightharpoonup \ h \quad \mbox{in}  \quad W^{1,0}_2(Q^+(\frac 9{10})),
\end{array}
\label{Weak_Convergence}
\end{equation}
\begin{equation}
\begin{array}c
u^m \ \to \ u \quad \mbox{in}  \quad L_3 (Q^+( \frac 9{10})), \\
h^m \ \to \ h \quad \mbox{in}  \quad L_3(Q^+( \frac 9{10})), \\
a^m\ \to \  a \quad \mbox{in}  \quad \mathbb R^3.\qquad \  \qquad
\end{array}
\label{Strong_Convergence}
\end{equation}

\medskip
\noindent
{\bf 5.} Passing to the limit in the equations (\ref{MHD_NS_m}), (\ref{MHD_Magnetic_m}), (\ref{BC_m}) we obtain
\begin{equation}
\gathered
\cd_t u   - \Delta u +  \nabla q = \rot h \times   a \qquad \mbox{in}\quad Q^+,
\\
 \div u =0 \qquad \mbox{in}\quad Q^+,
 \\
u|_{x_3=0}=0,
\endgathered
\label{LinProb1}
\end{equation}
\begin{equation}
\gathered
\cd_t h     -  \Delta h = \rot (u\times  a) \qquad \mbox{in}\quad Q^+,
\\
\div h =0 \qquad \mbox{in}\quad Q^+,
\\
\begin{array}c
h_3|_{x_3=0}=0, \qquad \frac{\cd h_1}{\cd x_3}\big|_{x_3=0}\ =\ \frac{\cd h_2}{\cd x_3}\big|_{x_3=0} \ = \ 0.
\end{array}
\endgathered
\label{LinProb2}
\end{equation}

\medskip
\noindent
{\bf 6.} From  (\ref{Strong_Convergence}) we conclude
$$
\lim\limits_{m\to +\infty} Y_\tau(u^m) = Y_\tau(u), \qquad \lim\limits_{m\to +\infty} \tilde Y_\tau(h^m) = \tilde Y_\tau(h)
$$
and hence
\begin{equation}
\limsup\limits_{m\to \infty}  Y_\tau (u^m, h^m , q^m ) \ \le \ Y_\tau(u) + \tilde Y_\tau(h) + \limsup\limits_{m\to \infty}  \hat Y_\tau (q^m).
\label{Y estimate}
\end{equation}
As $(u, h, q)$ is a solution to the linear problem (\ref{LinProb1}) --- (\ref{LinProb2}), from Theorem \ref{Linear_Estimate_Theorem} we obtain
\begin{equation}
Y_\tau(u) + \tilde Y_\tau(h) \ \le \ C(M)~ \tau^{1/3} ~Y_{1} (u, h, q)
\label{Y_u_h}
\end{equation}

\medskip
\noindent
{\bf 7.} Now we are going to estimate $\limsup\limits_{m\to \infty}  \hat Y_\tau (q^m)$. For this purpose we decompose $(u^m, q^m)$ and $(u,q)$ as
$$
\gathered
u^m = u^m_1 + u^m_2, \qquad q^m = q^m_1+q^m_2, \\
u=u_1+u_2,\ \qquad\quad q = q_1+q_2,
\endgathered
$$
where $(u^m_1, q^m_1)\in W^{2,1}_{\frac 98, \frac 32}(\Pi^+)\times W^{1,0}_{\frac 98, \frac 32}(\Pi^+)$  are determined as a solutions of the following initial boundary-value problems in $\Pi^+=\mathbb R^3_+\times (-1, 0) $:
$$
\gathered
\cd_t u^m_1 - \Delta u^m_1 + \nabla q_1^m = f^m \qquad \mbox{in}\quad \Pi_+, \\
\div u^m_1 =0 \qquad \mbox{in}\quad \Pi_+, \\
u^m_1|_{t=-1} =0, \qquad u^m_1|_{x_3=0}=0,
\endgathered
$$
where $f^m$ is defined by the expression $\rot h^m \times (\ep^m h^m+a^m)- \ep^m u^m \cdot \nabla u^m$ on the set $Q^+(\frac 9{10})$ and extended by zero onto the whole $\Pi^+$. Similarly,   $(u_1, q_1)$ are determined by the relations
\begin{equation}
\gathered
\cd_t u_1 - \Delta u_1 + \nabla q_1 = f \qquad \mbox{in}\quad \Pi_+, \\
\div u_1 =0 \qquad \mbox{in}\quad \Pi_+, \\
u_1|_{t=-1} =0, \qquad u_1|_{x_3=0}=0,
\endgathered
\label{Global_Linear_Problem}
\end{equation}
where $f$ determined by the expression $\rot h\times a $ on the set $Q^+(\frac 9{10})$ and extended by zero onto the whole $\Pi^+$.

\medskip
\noindent
{\bf 8.}  As functions $u^m_1-u_1$, $q^m_1-q_1$ are the solution of the first initial boundary-value problem in $\Pi^+$ with the right-hand side $f^m-f$ and zero initial and boundary conditions, we obtain the estimate (see, for example,   \cite{Solonnikov_Uspekhi}, Proposition 2.1)
\begin{equation}
\gathered
\| u^m_1 \|_{W^{2,1}_{\frac 98, \frac 32}(\Pi^+)} + \| \nabla q^m_1  \|_{L_{\frac 98, \frac 32}(\Pi^+) } \ \le C \| f^m\|_{L_{\frac 98, \frac 32}(Q^+(\frac 9{10}))}
\\
\| u^m_1-u_1 \|_{W^{2,1}_{\frac 98, \frac 32}(\Pi^+)} + \| \nabla q^m_1 - \nabla q_1 \|_{L_{\frac 98, \frac 32}(\Pi^+) } \ \le C \| f^m-f\|_{L_{\frac 98, \frac 32}(Q^+(\frac 9{10}))}
\endgathered
\label{u_1 convergence}
\end{equation}
Note that
\begin{equation}
\gathered
\| f^m\|_{L_{\frac 98, \frac 32}(Q^+(\frac 9{10}))} \ \le \ C(M)
\\
\| f^m - f\|_{L_{\frac 98, \frac 32}(Q^+(\frac 9{10}))} \ \to \ 0, \quad \mbox{as} \quad m\to \infty.
\endgathered
\label{f^m boundedness}
\end{equation}
So, taking into account the imbedding $W^{1,0}_{\frac 98, \frac 32}(Q^+(\frac 9{10})) \hookrightarrow L_{ \frac 32}(Q^+(\frac 9{10}))$  we can conclude that
$$
\begin{array}c
q^m_1\to q_1\quad \mbox{in} \quad L_{\frac 32}(Q^+(\frac 9{10}))
\end{array}
$$
and hence for any $\tau \in (0, \frac 9{10})$
$$
\lim\limits_{m\to \infty } Y_\tau (q^m_1 ) \ =  \ Y_\tau (q_1).
$$
On the other hand, $(u_1, q_1)$ is a solution of the linear Stokes problem in $Q^+$. Hence from Corollary \ref{Corollary_Linear_Estimate} we conclude
$$
Y_\tau (q_1) \ \le \  C(M)~\tau^{1/3}~ Y_{\frac 9{10}} (q_1)
$$

\medskip
\noindent
{\bf 9.} We need to estimate $Y_{\frac 9{10}} (q_1)$.
From imbedding theorem $L_{\frac 32}(B^+(\frac{9}{10})) \hookrightarrow W^1_{\frac 98}(B^+(\frac{9}{10}))$ we conclude
$$
Y_{\frac 9{10}} (q_1)\ \le \ C ~\| \nabla q_1\|_{L_{\frac 98,\frac 32}(B^+(\frac{9}{10}))}
$$
For the solution $(u_1, q_1)$ of the initial-boundary value problem (\ref{Global_Linear_Problem}) we have the estimate
$$
\| u_1\|_{W^{2,1}_{\frac 98, \frac 32}(Q^+(\frac 9{10}))} + \| \nabla q_1\|_{L_{\frac 98, \frac 32}(Q^+(\frac 9{10}))} \ \le \ C(M)~ \| \nabla h\|_{L_{\frac 98, \frac 32}(Q^+(\frac 9{10}))}
$$
Using H\" older inequality \ $\| \nabla h\|_{L_{\frac 98, \frac 32}(Q^+(\frac 9{10}))} \ \le \ C~\| \nabla h\|_{L_2(Q^+(\frac 9{10}))}$ and  taking into account the weak convergence (\ref{Weak_Convergence}) from which we conclude
$$
\| \nabla h\|_{L_2(Q^+(\frac 9{10}))} \ \le \ \liminf\limits_{m\to \infty} \| \nabla h^m \|_{L_2(Q^+(\frac 9{10}))},
$$
and using (\ref{Gradient boundedness}) we obtain
$$
Y_{\frac 9{10}} (q_1) \ \le \ C(M).
$$

\medskip
\noindent
{\bf 10.} Now we consider functions $(u^m_2, q^m_2)$  determined by relations
\begin{equation}
\gathered
u^m_2:=u^m-u^m_1, \qquad q^m_2:= q^m -q^m_1.
\endgathered
\label{u_2 definition}
\end{equation}
These functions satisfy the homogeneous Stokes problems in $Q^+(\frac{9}{10})$:
$$
\begin{array}c
\cd_t u^m_2 - \Delta u^m_2 + \nabla q_2^m = 0 \qquad \mbox{in}\quad Q^+(\frac{9}{10}), \\
\div u^m_2 =0 \qquad \mbox{in}\quad Q^+(\frac{9}{10}), \\
u^m_2|_{x_3=0}=0,
\end{array}
$$
$$
\begin{array}c
\cd_t u_2 - \Delta u_2 + \nabla q_2 = 0 \qquad \mbox{in}\quad Q^+(\frac{9}{10}), \\
\div u_2 =0 \qquad \mbox{in}\quad Q^+(\frac{9}{10}), \\
u_2|_{x_3=0}=0.
\end{array}
$$
Then
$$
\begin{array}c
\| u^m_2\|_{W^{2,1}_{9, \frac 32}(Q^+(\frac 4{5}))}+ \| \nabla q^m_2\|_{L_{9, \frac 32}(Q^+(\frac 4{5}))} \ \le \ C ~\Big(
\| u^m_2\|_{L_3(Q^+(\frac 9{10}))} + \| q^m_2\|_{L_{\frac 32}(Q^+(\frac 9{10}))}\Big)
\end{array}
$$
Note that due to (\ref{u_2 definition}), (\ref{Basic_boundedness}) and the first inequalities in (\ref{u_1 convergence}), (\ref{f^m boundedness}) we have the estimate
$$
\gathered
\| u^m_2\|_{L_3(Q^+(\frac 9{10}))} + \| q^m_2\|_{L_{\frac 32}(Q^+(\frac 9{10}))} \ \le \\
\le \
\| u^m\|_{L_3(Q^+(\frac 9{10}))} + \| q^m\|_{L_{\frac 32}(Q^+(\frac 9{10}))} \ +  \
\| u^m_1\|_{L_3(Q^+(\frac 9{10}))} + \| q^m_1\|_{L_{\frac 32}(Q^+(\frac 9{10}))}  \ \le  \\ \le \ C(M)
\endgathered
$$
On the other hand, from the H\" older inequality we obtain for any $\tau \in (0, \frac 45)$
$$
\gathered
\hat Y_{\tau}(q^m_2)= \tau \Big( \ - \!\!\!\!\!\!\!\!\! \int\limits_{Q^+(\tau)} |q_2^m- [q_2^m]_{B^+(\tau)}|^{\frac 32}~dx dt\Big)^{\frac 23} \ \le \
C \tau^2 ~ \Big( \ - \!\!\!\!\!\!\!\!\! \int\limits_{Q^+(\tau)} |\nabla q_2^m |^{\frac 32}~dx dt\Big)^{\frac 23} \\
\le \ C~\tau^{\frac 76} ~\| \nabla q^m_2 \|_{L_{9, \frac 32}(Q^+(\frac 45))} \ \le \ C(M)~\tau^{\frac 76}
\endgathered
$$

\medskip
\noindent
{\bf 11.} Summarizing all previous estimates we arrive at
$$
\gathered
\limsup\limits_{m\to \infty}  \hat Y_\tau (q^m) \ \le \  \lim\limits_{m\to \infty}  \hat Y_\tau (q^m_1) + \limsup\limits_{m\to \infty}  \hat Y_\tau (q^m_2) \
\le \ C(M) ~ \tau^{\frac 13}
\endgathered
$$
Taking into account (\ref{Y estimate})  and (\ref{Y_u_h}), we finally obtain
$$
\gathered
\limsup\limits_{m\to \infty}   Y_\tau (u^m, h^m, q^m) \ \le \   C(M) ~ \tau^{\frac 13}
\endgathered
$$
This estimate contradicts (\ref{Contradiction_asumptions}) whenever $C_* > C(M)$. Theorem \ref{Decay estimate theorem} is proved.

\bigskip
\noindent
Theorem \ref{CKN_theorem} follows from Theorem \ref{Decay estimate theorem}  in the standard way by iterations of the estimate (\ref{Decay_Estimate}), scaling arguments, and combination of boundary estimates with the internal estimates obtained in \cite{Vya}. See  details in  \cite{Seregin_JMFM}, \cite{Seregin_Aa}, \cite{SSS}, \cite{Seregin_Handbook}.

\newpage

\

\newpage


\section{Estimates of  Solutions \\ to  the Heat Equation}
\setcounter{equation}{0}

\bigskip
In this section we study solutions of the heat equations  with the lower order terms:
$$
\gathered
\cd_t H \ - \ \Delta H \ =  \ \div (v\otimes   H - H \otimes v) \qquad\mbox{in}\quad Q^+.
\\
v|_{x_3=0}=0,
\\
H_3|_{x_3=0}=0, \qquad  H_{\al, 3} |_{x_3=0}=0, \qquad \al=1,2.
\endgathered
$$
Namely, we assume the functions $(v, H)$ possess the following properties:
\begin{equation}
\gathered
v, \ H\in W^{1,0}_2(Q^+), \\ v|_{x_3=0}=0, \quad H_3|_{x_3=0}=0 \quad \mbox{in the sense of traces}, \\
\endgathered
\label{Class}
\end{equation}
 for any $\eta \in C^\infty_0(Q; \mathbb R^3)$ such that  $\eta_3|_{x_3=0}=0$ the following integral identity holds
\begin{equation}
\gathered
\int\limits_{Q^+} \Big( - H\cdot \cd_t \eta + \nabla H:\nabla \eta \Big)~dxdt \ =  \ - \int\limits_{Q^+} G:\nabla \eta~dxdt,  \\
\endgathered
\label{Integral_Identity}
\end{equation}
where $G\ = \ v\otimes H - H\otimes v$, and
\begin{equation}
\div v=0, \qquad\div H=0 \qquad\mbox{a.e. in}\quad Q^+.
\label{Divergent-free}
\end{equation}

\medskip
We start with the auxiliary results:
\begin{theorem} Assume conditions (\ref{Class}), (\ref{Integral_Identity}) hold. Denote by $v^*$ and $H^*$ the extensions of the functions $v$ and $H$ from $Q^+$ onto $Q$ obtained in the following way: the components $v^*_\al$, $H^*_\al$, $\al=1,2$ are the even extensions of the components $v_\al$, $H_\al$, $\al=1,2$, and the components $v^*_3$, $H^*_3$ are the odd extensions of $v_3$, $H_3$. Define also the function
$$
G^*= v^*\otimes   H^* - H^* \otimes v^*.
$$
Then the following relation holds:
\begin{equation}
\cd_t H^* - \Delta H^* \ = \  \div G^* \qquad\mbox{in}\qquad \mathcal D'(Q).
\label{g^*}
\end{equation}
Moreover, if (\ref{Divergent-free}) is satisfied, then
$$
\div v^* =0, \qquad \div H^* =0, \qquad \mbox{in}\quad \mathcal D'(Q).
$$
\label{Extension}
\end{theorem}

\noindent
{\bf Proof of Theorem \ref{Extension}: } The result is a direct consequence of the boundary conditions $v|_{x_3=0}=0$ and $H_3|_{x_3=0}=0$. Theorem \ref{Extension} is proved.

\medskip
Now we introduce the following functionals
\begin{equation}
\gathered
E(r) \ = \
\Big( ~\frac 1{r} \int\limits_{Q^+(r)} |\nabla v|^2~dxdt ~ \Big)^{1/2},
\\
 E_*(r) \ = \  \Big( \frac 1{r} \int\limits_{Q^+(r)} |\nabla H|^2~dxdt \Big)^{1/2},
\\
F_q( r) \ = \ \Big( \frac 1{r^{5-q}} \int\limits_{Q^+(r)} |H|^q~dxdt
\Big)^{1/q}.
\endgathered
\label{Definition E and F}
\end{equation}

\begin{theorem}
Assume conditions (\ref{Class}) --- (\ref{Divergent-free}) hold. Then for any $r\in (0, \frac 12)$ the following estimate holds
\begin{equation}
\gathered
\| H\|_{L_{1, \infty}(Q^+(r))} \ \le  \ c r^2~\Big(  1 + E(2r)\Big)\Big( F_2(2r)  + E_*(2r) + r \Big)
\endgathered
\label{H_{1,infty}}
\end{equation}
\label{L_{1,infty}}
\end{theorem}
\noindent
{\bf Proof of Theorem \ref{L_{1,infty}}: }
 Cosider the function
$$\eta=\frac{\zeta H^*}{(1+|H^*|^2)^{1/2}},$$ where $\zeta\in C_0^\infty( Q(2r))$  is a standard cut-off function.
The following relations are true:
$$
\gathered
\int\limits_{Q(2r)} \nabla H^* : \nabla \eta ~dxdt \ = \
\int\limits_{Q(2r)} \frac{\nabla H^* : H^* \otimes \nabla \zeta}{(1+|H^*|^2)^{\frac 12} } ~dxdt \ +  \\
+ \
\int\limits_{Q(2r)} \zeta \left(\frac{|\nabla H^*|^2}{ (1+|H^*|^2)^{\frac 12}} - \frac{|\nabla| H^*|^2|^2}{ 4(1+|H^*|^2)^{\frac 32}} \right)  ~dxdt ,
\endgathered
$$
and
$$
\begin{array}c
(1+|H^*|^2 ) |\nabla H^*|^2 - \frac 14 |\nabla| H^*|^2|^2 \ \ge \ |\nabla H^*|^2.
\end{array}
$$
Hence we obtain the estimate
$$
\gathered
\int\limits_{Q(2r)} \nabla H^* : \nabla \eta ~dxdt \ \ge \ \int\limits_{Q(2r)} \zeta \frac{|\nabla H^*|^2}{ (1+|H^*|^2)^{\frac 32}}~dxdt -
\int\limits_{Q(2r)} |\nabla H^* | |\nabla \zeta| ~dxdt .
\endgathered
$$
 Testing the equation (\ref{g^*}) by the function $\eta$ we obtain
$$
\gathered
\sup\limits_{t\in (-r^2,0)}\int\limits_{B(2r)} \zeta (1+ |H^*|^2)^{\frac 12}~dx + \int\limits_{Q(2r)} \zeta \frac{|\nabla H^*|^2}{ (1+|H^*|^2)^{\frac 32}}~dxdt   \le \\ \le \ c ~\int\limits_{Q(2r)} |\cd_t\zeta | (1+|H^*|^2)^{\frac 12}~dxdt \ + \ \int\limits_{Q(2r)} |\nabla H^* | |\nabla \zeta| ~dxdt \ + \\
+\ \int\limits_{Q(2r)}  \div G^*\cdot \eta~dxdt.
\endgathered
$$
Note that
$$
\div G^*  \ =  \  (H^*\cdot \nabla )v^* - (v^*\cdot \nabla ) H^*\qquad \mbox{a.e. in} \quad Q.
$$
Integrating by parts we obtain
$$
\gathered
\int\limits_{Q(2r)} (v^*\cdot \nabla ) H^*\cdot \eta ~dxdt \ = \
\int\limits_{Q(2r)} \zeta v^*\cdot \nabla (1+|H^*|^2)^{1/2}~dxdt \ = \\ = \ -
\int\limits_{Q(2r)}  v^*\cdot \nabla \zeta (1+|H^*|^2)^{1/2}~dxdt
\endgathered
$$
Hence
$$
\gathered
\int\limits_{Q(2r)} \div G^*\cdot \eta ~dxdt \ \le  \  c~ \| \nabla v\|_{L_2(Q^+(2r))} \| H\|_{L_2(Q^+(2r))}   \ +  \\
+ \ \frac {c}{r}~\| v\|_{L_2(Q^+(2r))} \| H \|_{L_2(Q^+(2r))} \ + \  \frac cr ~\| v\|_{L_1(Q^+(2r))}
\endgathered
$$
Taking into account  the Poincare inequality for $v$  and the H\" older inequality, we obtain  the estimate
$$
\gathered
\| H\|_{L_{1, \infty}(Q^+(r))} \ \le  \ c ~\Big(r^{\frac 12}  + \| \nabla v\|_{L_2(Q^+(2r))}\Big)\| H\|_{L_2(Q^+(2r))} \ + \\  +
\
c r^{\frac 32}\| \nabla H\|_{L_2(Q^+(2r))} \ +
\
c \Big(r^3  +  r^{\frac 52} \| \nabla v\|_{L_2(Q^+(2r))}\Big)
\endgathered
$$
which implies (\ref{H_{1,infty}}). Theorem \ref{L_{1,infty}} is proved.

\medskip
\noindent
The main result of this section is following:

\begin{theorem} Assume conditions (\ref{Class}) --- (\ref{Divergent-free}) hold.   Then there exist absolute positive constants $\ep_1$, $\al$ and $c$ such that for any $\ep\in (0,\ep_1)$ and any $K>0$ if
\begin{equation}
\sup\limits_{r\in (0,1)} E(r) < \ep \qquad\mbox{and}\qquad \sup\limits_{r\in (0,1)} E_*(r) < K
\label{Assumptions}
\end{equation}
then for any $0<r<\rho\le 1$
\begin{equation}
 F_2(r) \ \le \   c\left(\frac{r}{\rho}\right)^{\al} F_2(\rho ) \ + \   c \ep (K+1) .
\label{F_2}
\end{equation}
\label{Bound_F_2}
\end{theorem}

\noindent
{\bf Proof of Theorem \ref{Bound_F_2}: }

\bigskip
\noindent
{\bf 1.}
Denote by $v^*$ and $H^*$ the extensions of functions $v$ and $H$ from $Q^+$ onto $Q$ described in Theorem \ref{Extension}.
Fix  arbitrary $r\in (0, 1)$ and let $\zeta\in C^\infty(\bar Q)$ be a cut off function such that $\zeta \equiv 1$ on $Q(r)$ and $\operatorname{supp} \zeta\subset B\times (-1,0]$. Denote $ \Pi=\mathbb R^3\times (-1,0)$ and denote by ${\hat G}$ the function which coincides with   $G^*$ on $Q(\frac r2)$ and additionally possesses the following properties: ${\hat G}\in W^{1,0}_1(\Pi)\cap L_{\frac {18}{11}, \frac 65}(\Pi)$,
${\hat G}$ is compactly supported in $\Pi$,  and
\begin{equation}
\| {\hat G}\|_{L_{\frac {18}{11}, \frac 65}(\Pi)} \ \le \ c \| G^{*}\|_{L_{\frac {18}{11}, \frac 65}(Q(\frac r2))} \ \le \ c \| G\|_{L_{\frac {18}{11}, \frac 65}(Q^+(\frac r2))}
\label{G**}
\end{equation}

\medskip
\noindent
{\bf 2.}   We decompose $H^*$ as
$$
H^* \ =  \hat H \ + \tilde H,
$$
where $\hat H$ is a solution of the  Cauchy problem for the heat equation
\begin{equation}
\left\{ \ \begin{array}l
\cd _t \hat H -\Delta \hat H  \ = \ \div {\hat G} \qquad \mbox{in}\quad \Pi, \\
\hat H |_{t=-1} =0,
\end{array}\right.
\label{Cauchy problem}
\end{equation}
defined by the formula
$\hat H = \Ga * \div  {\hat G}= -\nabla \Ga *  {\hat G}$, where $\Ga$ is the fundamental solution of the heat operator.
The function $\tilde H$ satisfies the homogeneous heat equation
\begin{equation}
\begin{array}c
\cd_t \tilde H - \Delta \tilde H =  0 \qquad \mbox{in}\quad Q( \frac r2).
\end{array}
\label{Heat-2}
\end{equation}

\medskip
\noindent
{\bf 3.} Take arbitrary $\theta \in (0,\frac 12)$.
We estimate $\| H\|_{L_2(Q^+(\theta r))}$ in the following way
\begin{equation}
\gathered
\| H\|_{L_2(Q^+(\theta r))} \ \le \ \| H^* \|_{L_2(Q(\theta r))}  \ \le
\ \| \hat H\|_{L_2(Q(\theta r))} \ + \ \| \tilde H\|_{L_2(Q(\theta r))},
\endgathered
\label{F-2}
\end{equation}
For  $\| \hat H\|_{L_2(Q(\theta r))}$ we have
\begin{equation}
\| \hat H\|_{L_2(Q(\theta r))} \ \le \ c~\| \hat H\|_{L_2(Q(\frac  r2))}.
\label{F-3}
\end{equation}
As $\tilde H$ satisfies (\ref{Heat-2}) by local estimate of the maximum of $\tilde H$ via its $L_2$--norm  we obtain
\begin{equation}
\gathered
\| \tilde H\|_{L_2(Q(\theta r))} \ \le \ c ~\theta^{\frac {5}{2}} ~\| \tilde H\|_{L_2(Q(\frac r2))} \ \le
\\ \le
\ c ~\theta^{\frac {5}{2} } ~( \|  H^*\|_{L_2(Q(r))} + \| \hat H\|_{L_2(Q(\frac r2))})
\endgathered
\label{F-4}
\end{equation}

\medskip
\noindent
{\bf 4.} So, we need to estimate $\| \hat H\|_{L_2(Q(\frac r2))}$. As singular integrals are bounded on the anisotropic Lesbegue space $L_{s,l}$ (see, for example, \cite{Solonnikov_Uspekhi}) for the convolution  $\hat h = \Ga* {\hat G} $ we obtain the estimate
$$
\| \hat h \|_{W^{2,1}_{\frac{18}{11}, \frac 65}(Q(r))} \ \le  \ c \|  {\hat G} \|_{L_{\frac{18}{11}, \frac 65}(\Pi)}.
$$
On the other hand, from  the 3D-- parabolic imbedding theorem (see \cite{Besov})
$$
W^{2,1}_{s,l}(Q)\hookrightarrow W^{1,0}_{p,q}(Q), \quad\mbox{as}\quad 1-\left(\frac 3s+\frac 2l - \frac 3p -\frac 2q\right) \ge 0,
$$
for $p=q=2$ and $s=\frac {18}{11}$, $l=\frac 65$ and for $\hat H=- \nabla \hat h$ we obtain
$$
\| \hat H \|_{L_2(Q(r))} \ \le \ c ~\|  {\hat G} \|_{L_{\frac{18}{11}, \frac 65}(\Pi)}.
$$
(Note that the constant $c$ in this inequality does not depend on $r$). Taking into account (\ref{G**}) we arrive at
\begin{equation}
\| \hat H \|_{ L_2(Q(r))} \ \le \ c ~ \| G \|_{L_{\frac{18}{11}, \frac 65}(Q^+(\frac r2))}.
\label{imbedding}
\end{equation}

\medskip
\noindent
{\bf 5.}
From the definition of $G$ we obtain
$$
\gathered
\| G \|_{L_{\frac{18}{11}, \frac 65}(Q^+(\frac r2))} \ \le \ c ~
\left(\ \int\limits_{-r^2/4}^{0} \| v\otimes H\|_{L_{\frac {18}{11}}(B^+(r/2))}^{\frac 65}~dt\ \right)^{\frac 56}
\endgathered
$$
Applying the H\" older inequality and Sobolev imbedding $W^1_2(B^+(r))\hookrightarrow L_6(B^+(r))$ for $v$ we obtain
$$
\gathered
\| G \|_{L_{\frac{18}{11}, \frac 65}(Q^+(\frac r2))} \ \le \ c ~
\left(\ \int\limits_{-r^2/4}^{0} \| v\|_{L_6(B^+(r))}^{\frac 65} \| H\|_{L_{\frac {9}{4}}(B^+(r/2))}^{\frac 65}~dt\ \right)^{\frac 56} \ \le \\
\le \ c~ \left(\ \int\limits_{-r^2/4}^{0} \| \nabla v\|_{L_2(B^+(r))}^{\frac 65} \| H\|_{L_{\frac {9}{4}}(B^+(r/2))}^{\frac 65}~dt\ \right)^{\frac 56}.
\endgathered
$$
Interpolating $L_{\frac 94}$-- norm between $L_{1}$ and $L_6$ and  using the imbedding $W^1_2(B^+(r))\hookrightarrow L_6(B^+(r))$ again we obtain
$$
\gathered
\| H\|_{L_{\frac {9}{4}}(B^+(r/2))} \ \le \| H\|_{L_1(B^+(r/2))}^{\frac 13} \| H\|_{L_{6}(B^+(r/2))}^{\frac 23} \ \le \\ \le \
\| H\|_{L_1(B^+(r/2))}^{\frac 13} \Big(\| \nabla H\|_{L_{2}(B^+(r))}^{\frac 23} +  r^{-\frac 23}\| H\|_{L_{2}(B^+(r))}^{\frac 23}\Big)
\endgathered
$$
Hence we arrive at
$$
\gathered
\| G \|_{L_{\frac{18}{11}, \frac 65}(Q^+(\frac r2))} \
\le \ c~ \| H\|_{L_{1, \infty}(Q^+(r/2))}^{\frac 13} \ \times  \\ \times \ \left(\ \int\limits_{-r^2/4}^{0} \| \nabla v\|_{L_2(B^+(r))}^{\frac 65} \Big(\| \nabla H\|_{L_{2}(B^+(r))}^{\frac 45} + r^{-\frac 45} \|H\|_{L_{2}(B^+(r))}^{\frac 45}\Big)
~dt\ \right)^{\frac 56}.
\endgathered
$$
Applying the H\" older inequality we obtain
$$
\gathered
\| G \|_{L_{\frac{18}{11}, \frac 65}(Q^+(\frac r2))} \
\le \ c~ \| H\|_{L_{1, \infty}(Q^+(r/2))}^{\frac 13}
\| \nabla v \|_{L_2(Q^+(r))} \ \times \\ \times \   \Big( \| \nabla H \|_{L_2(Q^+(r))}^{\frac 23} + r^{-\frac 23} \| H \|_{L_2(Q^+(r))}^{\frac 23}\Big)
\endgathered
$$

\medskip
\noindent
{\bf 6.} Estimating $\| H\|_{L_{1, \infty}(Q^+(r/2))}$ using   Theorem \ref{L_{1,infty}} we obtain
\begin{equation}
\gathered
\| G \|_{L_{\frac{18}{11}, \frac 65}(Q^+(\frac r2))} \
\le \
c r^{\frac 23}~\Big(  1 + E(r)\Big)^{\frac 13}\Big( F_2(r) + E_*(r) +  r \Big)^{\frac 13}
\ \times \\ \times \  \| \nabla v \|_{L_2(Q^+(r))}  \Big( \| \nabla H \|_{L_2(Q^+(r))}^{\frac 23} + r^{-\frac 23} \| H \|_{L_2(Q^+(r))}^{\frac 23}\Big)
\endgathered
\label{Estimate G}
\end{equation}

\medskip
\noindent
{\bf 7.}
Gathering estimates (\ref{F-2}) --- (\ref{Estimate G}) together we arrive at
$$
\gathered
\| H \|_{L_2(Q^+(\theta r))} \ \le \  c ~ \theta^{\frac 52}  \| H \|_{L_2(Q^+( r))}  \ +  \
c ~  r^{\frac 23} \| \nabla v \|_{L_2(Q^+(r))}   \ \times \\ \times \
\Big(  1 + E(r)\Big)^{\frac 13}\Big( F_2(r)  +  E_*(r) +  r \Big)^{\frac 13}  \Big( \| \nabla H \|_{L_2(Q^+(r))}^{\frac 23} + r^{-\frac 23} \| H \|_{L_2(Q^+(r))}^{\frac 23}\Big)
\endgathered
$$
Dividing this inequality by $(\theta r)^{\frac 32}$ we arrive at
$$
\gathered
F_2(\theta r) \ \le \  c~\theta ~ F_2(r) \ + \\
+ \ c(\theta) ~E(r) \Big(  1 + E(r)\Big)^{\frac 13}\Big( F_2(r)  +  E_*(r) + r \Big).
\endgathered
$$
Hence we obtain
$$
\gathered
F_2(\theta r) \ \le \ \Big( ~c\theta + c(\theta) E(r) \big(  1 + E(r)\big)^{\frac 13} ~\Big) F_2(r) \ + \\
+ \ c(\theta) ~E(r) \Big(  1 + E(r)\Big)^{\frac 13}  \left( E_*(r) + 1 \right).
\endgathered
$$
Taking into account assumptions (\ref{Assumptions}) for $\ep<1$ we obtain
$$
\gathered
F_2(\theta r) \ \le \ \Big( ~ c\theta + c(\theta) \ep ~\Big) F_2(r) \
+ \ c(\theta) \ep   ( K +  1),
\endgathered
$$
valid for any $r\in (0,1)$ and any $\theta\in (0,\frac 12]$.

\medskip
\noindent
{\bf 8.} Choosing $\theta\in (0,\frac 12]$ so that
$$
 c\theta  \ = \ \frac 14
$$
and then choosing $\ep_1\in (0,1)$ so that
$$
\frac 14 + c(\theta)\ep_1 \ \le \ \frac 12
$$
we obtain the estimate
$$
F_2(\theta r) \ \le \ \frac 12 ~F_2(r)  \ + \ c \ep   ( K +  1).
$$
Iterating this estimate we derive (\ref{F_2}). Theorem \ref{Bound_F_2} is proved.

\newpage

\section{Estimates of Energy Functionals}
\setcounter{equation}{0}

\bigskip

In the previous section we defined functionals  $F_q(r)$,  $E(r)$, and $E_*(r)$, see (\ref{Definition E and F}).
Now we define few more functionals. Note that all these functionals are invariant with respect to the natural scaling of the MHD system. For $r\le 1$, $q\in [1,\frac{10}3]$ and $s\in [1,\frac 98]$ we introduce the following quantities:
$$
\begin{array}c
A( r) \equiv \Big( \frac 1{r}
\sup\limits_{t\in (-r^2, 0)} \int\limits_{B^+(r)} |v|^{2}~dy
\Big)^{1/2}, \\
A_*( r) \equiv \Big( \frac 1{r}
\sup\limits_{t\in (-r^2, 0)} \int\limits_{B^+(r)} |H |^{2}~dy
\Big)^{1/2},
\\
C_q( r) \equiv \Big( \frac 1{r^{5-q}} \int\limits_{Q^+(r)} |v|^q~dydt
\Big)^{1/q},
\\ D(r) \equiv \Big( \frac 1{r^2} \int\limits_{Q^+(r)}
|p - [p]_{B^+(r)}|^{3/2}~dydt \Big)^{2/3},
\\
D_s(r) = R^{\frac 53 -\frac 3s} \Big( \int\limits_{-r^2}^{0} \Big(
\int\limits_{B^+(r)}  |\nabla p|^{s}~dy \Big)^{\frac 1s \cdot
\frac 32}~dt \Big)^{2/3},
\end{array}
$$
$$
C(r) = C_3(r), \qquad F(r)= F_3(r), \qquad D_*(r)= D_{\frac{36}{35}}(r).
$$

\noindent
First we formulate the set of results following from the general theory of functions:

\begin{theorem} \label{Interpolation}
Assume $v$, $H\in W^{1,0}_2(Q^+)$ and $p\in W^{1,0}_{\frac 98, \frac 32}(Q^+)$ are arbitrary functions. Assume $v|_{x_3=0}=0$. Then the following inequalities hold:
\begin{equation}\label{C_3}
C(r) \ \le  \ ~A^{\frac 12}(r)E^{\frac 12}(r),
\qquad
F(r)\ \le \ A_*^{\frac 12}(r)[ E_*^{\frac 12}(r) +  F_2^{\frac 12}(r)]
\end{equation}
\begin{equation}
D(r)\ \le \  c D_1(r), \qquad D_1(r)\ \le \ c D_s(r), \qquad \forall s>1.
\label{D}
\end{equation}

\end{theorem}

\noindent
{\bf Proof of Theorem \ref{Interpolation}:} The proof follows from interpolation inequalities and imbedding theorems.
Proof of the similar inequalities for the Navier-Stokes system can be found in \cite{LS}.

\bigskip\noindent
Now we formulate a theorem concerning boundary suitable weak solutions to the MHD system.

\begin{theorem}\label{Estimate_F_2}
Assume $(v,H, p)$ is a boundary suitable weak solution to the MHD equations in $Q^+$.
Then for any $r\in (0,1)$ and $\theta \in (0,\frac 12)$ the following inequalities hold
\begin{equation}
\gathered
A(r/2) + A_*(r/2) + E(r/2) + E_*(r/2) \ \le \\ \le \ c~ \Big( C_2(r)+ F_2(r) + C^{\frac 12}(r)D^{\frac 12}(r) +  C^{\frac 32}(r)\Big) \ + \\
+ \ c~ \Big( C^{\frac 12}(r)A_*^{\frac 12}(r)E_*^{\frac 12}(r) + F^{\frac 12}(r)A_*^{\frac 12}(r)E^{\frac 12}(r) \Big)
\label{Lo_En_Iq}
\endgathered
\end{equation}
\begin{equation}
\gathered
D_*(\theta r) \ \le \ c ~\theta^{\frac 43} ~\Big( D_*(r) + E(r)  \Big) \ + \\ +  \ c(\theta) ~\Big(  A^{\frac 23}(r) E^{\frac 43}(r) +  A^{\frac 56}_*(r)  F^{\frac 16}(r) E_*(r) \Big)
\endgathered
\label{D_*}
\end{equation}

\end{theorem}

\bigskip
\noindent
{\bf Proof of Theorem \ref{Estimate_F_2}, estimate (\ref{Lo_En_Iq}):}

\bigskip
\noindent
{\bf 1.} Estimate (\ref{Lo_En_Iq}) follows from (\ref{LEI}) in a standard way. We just explain the specific estimates of the terms
$$
I_1:=\int\limits_{Q^+(r)} ~ |H|^2 (v \cdot \nabla \zeta)   ~dxdt  \quad \mbox{and}   \quad I_2:=\int\limits_{Q^+(r)} ~ (v \cdot H) (H\cdot \nabla \zeta)  ~dxdt.
$$

\medskip
\noindent
{\bf 2.} $I_1$ we transform in the following way
$$
\begin{array}c
I_1\ = \ \int\limits_{Q^+(r)} ~ \Big( |H|^2  - [|H|^2]_{B^+(r)} \Big) (v \cdot \nabla \zeta)  ~dxdt
\end{array}
$$
Applying the H\" older  inequality we obtain
$$
|I_1| \ \le \ \frac cr ~\int\limits_{-r^2}^0 \left\| |H|^2  - [|H|^2]_{B^+(r)}\right\|_{L_{\frac 32}(B^+(r))} \| v\|_{L_{3}(B^+(r))}~dt
$$
Applying the inequality $\| f-[f]_{B^+(r)}\|_{L_{\frac 32}(B^+(r))} \le c \| \nabla f \|_{L_1(B^+(r))}$, we arrive at
$$
\gathered
|I_1| \ \le \ \frac cr ~\int\limits_{-r^2}^0 \| \nabla |H|^2  \|_{L_{1}(B^+(r))} \| v\|_{L_{3}(B^+(r))}~dt  \ \le \\ \le \
\frac cr ~\int\limits_{-r^2}^0 \| H\|_{L_2(B^+(r))} \| \nabla H  \|_{L_{2}(B^+(r))} \| v\|_{L_{3}(B^+(r))}~dt  \ \le \\ \le \
\frac c{r^{2/3}}  ~\| H\|_{L_{2,\infty}(Q^+(r))} \| \nabla H \|_{L_2(Q^+(R))} \| v \|_{L_3(Q^+(r))}\ \le \  cr~ A_*(r) E_*(r) C(r)
\endgathered
$$

\medskip
\noindent
{\bf 3.} For $I_2$ we obtain relations
$$
I_2 \ = \ \int\limits_{Q^+(r)} ~ \Big ( (v \cdot H)- [v \cdot H]_{B^+(r)}\Big) (H\cdot \nabla \zeta)  ~dxdt
$$
Hence
$$
\gathered
|I_2| \ \le   \ \frac cr ~\int\limits_{-r^2}^0 \left\| (v \cdot H)- [v \cdot H]_{B^+(r)} \right\|_{L_{2}(B^+(r))} \| H\|_{L_{2}(B^+(r))}~dt \ \le \\
\le \ \frac cr ~ \| H\|_{L_{2,\infty}(Q^+(r))} \int\limits_{-r^2}^0 \left\| \nabla (v \cdot H) \right\|_{L_{\frac 65}(B^+(r))} ~dt \ \le \
\frac cr ~ \| H\|_{L_{2,\infty}(Q^+(r))} \ \times \\ \times \ \int\limits_{-r^2}^0 \Big( \| \nabla v \|_{L_2(B^+(r))} \| H \|_{L_3(B^+(r))} + \| \nabla H \|_{L_2(B^+(r))} \| v \|_{L_3(B^+(r))}\Big)~dt \ \le \\
\le \ \frac c{r^{2/3}}  ~\| H\|_{L_{2,\infty}(Q^+(r))}  \Big( \| \nabla v \|_{L_2(Q^+(r))} \| H \|_{L_3(Q^+(r))} + \| \nabla H \|_{L_2(Q^+(r))} \| v \|_{L_3(Q^+(r))}\Big) \endgathered
$$
So, we obtain
$$
|I_2| \ \le \  cr~ A_*(r) ~\Big(~E(r) F(r) + E_*(r) C(r)~\Big)
$$

\bigskip
\noindent
{\bf Proof of Theorem \ref{Estimate_F_2}, estimate (\ref{D_*}):}

\medskip
\noindent
{\bf 1.} To obtain (\ref{D_*}) we apply the method developed in \cite{Seregin_JMFM}, \cite{Seregin_ZNS271}, see also \cite{SSS}. Denote $\Pi_r=\mathbb R^3_+\times (-r^2,0)$. We fix $r\in (0,1]$ and $\theta\in (0,\frac 12)$ and define a function $g:\Pi_{r}^+\to \mathbb R^3$ by the formula
$$
g \ = \  \left\{\begin{array}{cl}
 \rot H \times H-(v\cdot \nabla )v , & \mbox{in } \  Q^+(r), \\
0, & \mbox{in } \ \Pi_r^+\setminus Q^+(r)
\end{array}
\right.
$$
Then we decompose $v$ and $p$ as
$$
v\ = \ \hat v + \tilde v,\qquad p \ = \ \hat p +\tilde p,
$$
where $(\hat v, \hat p)$ is a solution of the Stokes initial boundary value problem in a half-space
$$
\gathered
\left\{\begin{array}c  \cd_t \hat v - \Delta \hat v+\nabla \hat p \ = \ g, \\ \div \hat v =0 \end{array}\right. \qquad \mbox{in}\quad \Pi_r^+, \\
\hat v|_{t=0}=0, \qquad \hat v|_{x_3=0}=0,
\endgathered
$$
and $(\tilde v, \tilde p)$ is a solution of the homogeneous Stokes system in $Q^+(r)$:
$$
\gathered
\left\{\begin{array}c  \cd_t \tilde v - \Delta \tilde v +\nabla \tilde p\ = \ 0, \\ \div \tilde v =0 \end{array}\right. \qquad \mbox{in}\quad Q^+(r),  \\
\tilde v|_{x_3=0}=0.
\endgathered
$$

\medskip
\noindent
{\bf 2.} For $\nabla \hat p$ and $\nabla \tilde p$ the following estimates hold (see \cite{Seregin_ZNS271}, see also \cite{Solonnikov_ZNS288}):
$$
\gathered
\| \nabla \hat p\|_{L_{\frac{36}{35},\frac 32}(Q^+(r))} \ + \ \frac 1r \| \nabla \hat v \|_{L_{\frac{36}{35},\frac 32}(Q^+(r))}   \ \le \\ \le \  c~\Big( ~\| H\times \rot H\|_{L_{\frac{36}{35}, \frac 32}(Q^+(r))} \ + \ \| (v\cdot \nabla)v \|_{L_{\frac{36}{35}, \frac 32}(Q^+(r))}~\Big),
\endgathered
$$
$$
\| \nabla \tilde p\|_{L_{\frac{36}{35},\frac 32}(Q^+(\theta r))} \ \le \ c~\theta^{\frac{31}{12}} ~\Big( ~\frac 1r \| \nabla \tilde v \|_{L_{\frac{36}{35},\frac 32}(Q^+(r))}
\ + \  \|  \nabla \tilde p \|_{L_{\frac{36}{35},\frac 32}(Q^+(r))}~\Big).
$$

\medskip
\noindent
{\bf 3.} From the H\" older inequality we obtain
$$
\gathered
\| H\times \rot H\|_{L_{\frac{36}{35}, \frac 32}(Q^+(r))} \ \le \ c ~r^{\frac 29} ~\| H \|_{L_{2,\infty}(Q^+(r))}^{\frac 56} \| \nabla H\|_{L_2(Q^+(r))} \|  H \|_{L_3(Q^+(r))}^{\frac 16}
\endgathered
$$
$$
\gathered
\| (v\cdot \nabla)v \|_{L_{\frac{36}{35}, \frac 32}(Q^+(r))} \ \le \ c~r^{\frac 14}~ \| (v\cdot \nabla )v \|_{L_{\frac{9}{8}, \frac 32}(Q^+(r))} \ \le  \\ \le \ c~r^{\frac 14}~
\| v\|_{L_{2,\infty}(Q^+(r))}^{\frac 23} \| \nabla v\|_{L_2(Q^+(r))}^{\frac 43}
\endgathered
$$
Representing $\tilde v = v-\hat v$, $\tilde p=p-\hat p$ and gathering all  above estimates for $\hat p$ and $\hat v$ we obtain
$$
\gathered
D_*( \theta r) \ \le \ c~\theta^{\frac 43} ~\Big( ~D_*(r) + E(r) + A^{\frac 23}(r) E^{\frac 43}(r) + A_*^{\frac 56}(r) E_*(r) F^{\frac 16}(r) ~\Big) \ +  \\
+ \ c(\theta)~ \Big( ~A^{\frac 23}(r) E^{\frac 43}(r) + A_*^{\frac 56}(r) E_*(r) F^{\frac 16}(r) ~\Big)
\endgathered
$$
Theorem \ref{Estimate_F_2} is proved.

\newpage

\section{CKN condition and \\ Partial Regularity of Solutions}
\setcounter{equation}{0}

\bigskip
In this section we present the proofs of Theorems \ref{CKN_theorem} and \ref{Partial_Regularity}.

\begin{theorem} \label{Boundedness_A_A_D}
Denote by \ $\mathcal E(r)$ \ the following functional
$$
\mathcal E (r) \ = \ A(r)+   A_*(r) +  D_*(r),
$$
and let $\ep_1>0$ be the absolute  constant defined in Theorem \ref{Bound_F_2}. For any $K>0$ there exists a constant $c(K)>0$ such that  for any $\ep\in (0,\ep_1]$ and any boundary suitable weak solution  $(v,H,p)$ of the MHD system in $Q^+$ if
\begin{equation}
\sup\limits_{r\in (0,1)}  E(r)\le  \ep, \qquad  \sup\limits_{r\in (0,1)} E_*(r)  \ \le \ K,
\label{Statement_of_Theorem_E_E1}
\end{equation}
and
\begin{equation}
F_2(1) \ \le \ M,
\label{M}
\end{equation}
then for any $0< r<\rho \le 1$
\begin{equation}
\mathcal E(r) \ \le \ c \left(\frac r\rho\right)^\be\mathcal E(\rho) \  +  \ c(K)(1+\rho^\al M).
\label{Estimate_varE}
\end{equation}
where $\be>0$ is some absolute constant.
%
\end{theorem}

\bigskip
\noindent
{\bf Proof of Theorem \ref{Boundedness_A_A_D}:}

\medskip\noindent
{\bf 1.} Without loss of generality we can assume $K\ge 1$. Then from  (\ref{F_2})  we obtain
$$
F_2(r)\le cr^\al M + c  K.
$$
 From this inequality and (\ref{C_3}) we obtain
\begin{equation}
C(r) \ \le \ c~ \mathcal E^{\frac 12} (r) \ep^{\frac 12}_1, \qquad F(r)\ \le \ c~ \mathcal E^{\frac 12} (r)\Big( K^{\frac 12} +  r^{\frac \al 2}M^{\frac 12}  \Big)
\label{C}
\end{equation}

\medskip\noindent
{\bf 2.} Assume $r\in (0,1)$ and $\theta \in (0,\frac 12)$.
From (\ref{Lo_En_Iq}) with the help of (\ref{D}) and the Young inequality  we obtain
$$
\gathered
\mathcal E(\theta r) \ \le \  c~\Big( F_2(2\theta r) + D_*(2\theta r)\Big) \  +
\\ + \
c(\theta)\Big( C_2(r) + C(r) + C^{\frac 32}(r) + C^{\frac 12}(r)A_*^{\frac 12}(r)E_*^{\frac 12}(r) + F^{\frac 12}(r)A_*^{\frac 12}(r)E^{\frac 12}(r)\Big)
\endgathered
$$
Taking into account (\ref{C}) and (\ref{Statement_of_Theorem_E_E1}) we obtain
\begin{equation}
\gathered
\mathcal E(\theta r) \ \le \  c~\Big( F_2(2\theta r) + D_*(2\theta r)\Big) \  +
\\ + \
c(\theta)\Big( \ep_1 + \mathcal E^{\frac 12} (r) \ep_1^{\frac 12} + \mathcal E^{\frac 34} (r) \ep_1^{\frac 34} + \ep_1^{\frac 14}\mathcal E^{\frac 34}(r)K^{\frac 12} +
 ( K^{\frac 14} +  r^{\frac \al 4}M^{\frac 14})
 \mathcal E^{\frac 34}(r)\ep_1^{\frac 12}\Big)
\endgathered
\label{Similar to this}
\end{equation}
Applying the Young inequality $ab\le \ep a^p+C_\ep b^{p'}$ we obtain
$$
\gathered
\mathcal E(\theta r)  \le   \frac 14 \mathcal E(r)  +   c\Big( F_2(2\theta r) + D_*(2\theta r)\Big)   +  c(\theta )c(K)  +   c(\theta) r^\al M.
\endgathered
$$

\medskip\noindent
{\bf 3.} From (\ref{F_2}) and (\ref{D_*}) we obtain
$$
\gathered
F_2(2\theta r) + D_*(2\theta r)  \le  c \theta^\al \Big(F_2(r) + D_*(r) \Big) +  c\ep_1(1+K) +  \\ +
c(\theta) ~\Big(  A^{\frac 23}(r) E^{\frac 43}(r) +  A^{\frac 56}_*(r)  F^{\frac 16}(r) E_*(r) \Big)
\endgathered
$$
Taking into account  (\ref{C}) and the obvious inequality $F_2(r)\le A_*(r)$ we arrive at
$$
\gathered
F_2(2\theta r) + D_*(2\theta r)  \ \le \  c \theta^\al \mathcal E(r) +  c(K) + \\ + c(\theta) \Big(\mathcal E^{\frac 23}(r) \ep_1^{\frac 43} +
\mathcal E^{\frac {11}{12}}(r) ( K^{\frac 1{12}} + r^{\frac \al{12}} M^{\frac 1{12}}) K\Big)
\endgathered
$$
Applying the Young inequality we get
$$
\gathered
F_2(2\theta r) + D_*(2\theta r)  \ \le \  \Big(\frac 14 + c \theta^\al\Big) \mathcal E(r) +  c(\theta)c(K) (1+
r^\al M)
\endgathered
$$

\medskip\noindent
{\bf 4.} Gathering the estimates we obtain
$$
\gathered
\mathcal E(\theta r)  \le    \Big(\frac 14 + c \theta^\al\Big) \mathcal E(r) +  c(\theta)c(K) (1+
r^\al M).
\endgathered
$$
Fixing $\theta\in (0,\frac 12)$ so that
$$
\frac 14 + c \theta^\al = \frac 12
$$
Hence
$$
\gathered
\mathcal E(\theta r)  \le    \frac 12\mathcal E(r) +  c(\theta)c(K) (1+
r^\al M).
\endgathered
$$
Iterating this inequality we obtain (\ref{Estimate_varE}). Theorem \ref{Boundedness_A_A_D} is proved.

\begin{theorem} \label{Estimate_A_A_D} Assume all conditions of Theorem \ref{Boundedness_A_A_D} hold and fix $\rho_0\in (0,1)$ so that
\begin{equation}
\rho_0^\al M \ \le \ 1.
\label{rho_0}
\end{equation}
Then for any $0<r<\rho\le \rho_0 $  the following estimates hold:
\begin{equation}
\gathered
 A(r) + A_*( r)    \le  c\left(\frac r\rho \right)^{\ga}
\Big( A( \rho ) + A_*( \rho) \Big)   +   \ep^{\frac 14}D(\rho) +  G(K, \ep)
\endgathered
\label{Statement_of_Theorem_A_A}
\end{equation}
\begin{equation}
\gathered
D(r) \ \le \ c\left( \frac{r}{\rho}\right)^{\ga} D(\rho)\ + \ c(K) \Big(A^{\frac {11}{12}}(\rho)+A_*^{\frac {11}{12}}(\rho)\Big) + G(K, \ep)
\endgathered
\label{Cor}
\end{equation}
where $\ga>0$ is some absolute constant and $G$ is a continuous function possessing the following property:
\begin{equation}
\gathered
\mbox{for any fixed}\quad K >0 \quad G(K, \ep)\to 0 \quad \mbox{as}\quad\ep\to 0.
\endgathered
\label{Property of G}
\end{equation}
\end{theorem}

\bigskip
\noindent
{\bf Proof of Theorem \ref{Estimate_A_A_D}:}

\medskip\noindent
{\bf 1.} From (\ref{C_3}) taking into account (\ref{rho_0}) we obtain
\begin{equation}
C(r) \ \le A^{\frac 12}(r)\ep^{\frac 12}, \qquad F(r) \ \le \  A_*^{\frac 12}(r)(K^{\frac 12}+1)
\label{F_K}
\end{equation}

\medskip\noindent
{\bf 2.} Take arbitrary $r\in (0,\rho_0)$ and $\theta\in (0,\frac 12)$.
Denote by \ $\mathcal E(r)$ \ the following functional
$$
\mathcal E_* (r) \ = \ A(r)+   A_*(r),
$$
Then from  (\ref{Lo_En_Iq}) similar to (\ref{Similar to this}) using (\ref{F_K}) we derive
$$
\gathered
\mathcal E_*(\theta r)  \ \le \  F_2(2\theta r) \ + \  C^{\frac 12}(2\theta r)D^{\frac 12}(2\theta r) \ +
\\
+\ c(\theta)\Big( \mathcal E_*^{\frac 12}(r)\ep^{\frac 12} + \mathcal E_*^{\frac 34}(r)\ep^{\frac 34} +
\mathcal E_*^{\frac 34}(r) K^{\frac 12}\ep^{\frac 14} + \mathcal E_*^{\frac 34}(r) (K^{\frac 14} +1) \ep^{\frac 12}\Big)
\endgathered
$$
Applying the Young inequality and using (\ref{D}) we obtain
\begin{equation}
\gathered
\mathcal E_*(\theta r)  \ \le \  \frac 18 ~\mathcal E_*(r) \ + \  c(\theta)G(K,\ep)  \ +
\\
+ \ F_2(2\theta r) \ + \  C^{\frac 12}(2\theta r)D^{\frac 12}_*(2\theta r)
\endgathered
\label{1}
\end{equation}

\medskip\noindent
{\bf 3.} From (\ref{F_2}) we conclude
\begin{equation}
F(2\theta r) \ \le c\theta^\al  \mathcal E_*(r) \ + \ G(K, \ep).
\label{2}
\end{equation}

\medskip\noindent
{\bf 4.}
From (\ref{D_*}) for $r\le \rho_0$ with the help of (\ref{F_K}) and the Young inequality we obtain
\begin{equation}
\gathered
D_*(2\theta r) \ \le \ c ~\theta^{\be} D_*(r) \ + \ c(\theta)c(K)\mathcal E^{\frac {11}{12}}_*(r)  + c(\theta)G(K,\ep)
\endgathered
\label{DDD}
\end{equation}
Hence from (\ref{F_K}) we obtain
$$
\gathered
 C^{\frac 12}(2\theta r)D^{\frac 12}(2\theta r) \ \le \ c(\theta) \mathcal E_*^{\frac 14}(r)\ep^{\frac 14} D_*^{\frac 12}(r)    \ +  \\
 + \ c(\theta) c(K) \ep^{\frac 14} \mathcal E^{\frac {17}{24}}_*(r) + c(\theta)G(K,\ep)
\endgathered
$$
Applying the Young inequality we arrive at
\begin{equation}
\gathered
 C^{\frac 12}(2\theta r)D^{\frac 12}(2\theta r) \ \le \ \frac 18 \mathcal E_*(r)   +  \frac 12 \ep^{\frac 14} D_*(r)  + \ c(\theta)G(K, \ep)
\endgathered
\label{3}
\end{equation}

\medskip\noindent
{\bf 5.} Gathering estimates (\ref{1}) --- (\ref{3}) we obtain the inequality
$$
\gathered
\mathcal E_*(\theta r)  \ \le \  \Big(\frac 14 + c\theta^\ga \Big) ~\mathcal E_*(r)  + \frac 12 \ep^{\frac 14} D_*(r) +  c(\theta)G(K,\ep)
\endgathered
$$
Choosing $\theta\in (0,\frac 12)$ so that
 $$
 \frac 14 + c\theta^\al \ = \  \frac 12
 $$
we obtain
$$
\gathered
\mathcal E_*(\theta r)  \ \le \  \frac 12 ~\mathcal E_*(r)  + \frac 12 \ep^{\frac 14} D_*(r) +  c(\theta)G(K,\ep)
\endgathered
$$
Iterating this inequality we obtain (\ref{Statement_of_Theorem_A_A}).

\medskip\noindent
{\bf 6.} Choosing in (\ref{DDD}) $\theta\in (0,\frac 12)$ so that
$$
c\theta^\be \ =  \ \frac 12
$$
and iteration the inequality we obtain we derive (\ref{Cor}). Theorem \ref{Estimate_A_A_D} is proved.

\begin{theorem}  \label{CKN_condition} For any $K>0$ there exists a constant $\ep_0(K)>0$ such that if the condition (\ref{Statement_of_Theorem_E_E1}) holds with $\ep\le \ep_0$,  then there exists $\rho_*\in (0,1)$ such that
$$
\begin{array}c
 \Big( C(\rho_*) + F(\rho_*) + D(\rho_*)\Big) \ < \ \ep_*^{\frac 13},
\end{array}
$$
where the constant $\ep_*>0$ is defined in Theorem \ref{Fixed_r}.
\end{theorem}

\bigskip
\noindent
{\bf Proof of Theorem \ref{CKN_condition}:}

\medskip\noindent
{\bf 1.} From (\ref{Estimate_varE}) we obtain
$$
\limsup\limits_{r\to 0} D_*(r) \le c(K).
$$

\medskip\noindent
{\bf 2.} From (\ref{Statement_of_Theorem_A_A}) we derive
$$
\gathered
\limsup\limits_{r\to 0} \Big(A(r)+A_*(r)\Big) \ \le  \ \ep^{\frac 14} \limsup\limits_{\rho\to 0} D(\rho) + G(K, \ep ) \ \le \\ \le  \ \ep^{\frac 14} c(K) + G(K, \ep).
\endgathered
$$

\medskip\noindent
{\bf 3.} From (\ref{Cor}) we obtain
$$
\gathered
\limsup\limits_{r\to 0} D_*(r) \ \le  \ c(K)\limsup\limits_{\rho\to 0} \Big(A^{\frac {11}{12}}(\rho)+ A^{\frac {11}{12}}_*(\rho)\Big)  + G(K, \ep ) \ \le \\ \le  \
c(K) \Big(\ep^{\frac 14} c(K) + G(K, \ep)\Big)^{\frac {11}{12}}
 + G(K, \ep).
\endgathered
$$

\medskip\noindent
{\bf 4.} From (\ref{C_3}) we conclude
$$
\gathered
\limsup\limits_{r\to 0}  \Big( C(r)+F(r)\Big) \ \le (\ep^{\frac 12} + K^{\frac 12}) \limsup\limits_{r\to 0} \Big(A(r)+A_*(r)\Big)  \ \le \\ \le \
(\ep^{\frac 12} + K^{\frac 12})\Big( \ep^{\frac 14} c(K) + G(K, \ep)\Big)^{\frac 12}.
\endgathered
$$

\medskip\noindent
{\bf 5.} Taking into account  (\ref{Property of G}) for any $K>0$ we can find $\ep_0(K)>0$ such that for any $\ep\in (0,\ep_0)$
$$
c(K) \Big(\ep^{\frac 14} c(K) + G(K, \ep)\Big)^{\frac {11}{12}}
 + G(K, \ep) \ < \ \frac{\ep_*^{\frac 13}}{2}
$$
and
$$
(\ep^{\frac 12} + K^{\frac 12})\Big( \ep^{\frac 14} c(K) + G(K, \ep)\Big)^{\frac 12} \ < \ \frac{\ep_*^{\frac 13}}{2}.
$$
Hence for $\ep\in (0,\ep_0)$
$$
\limsup\limits_{r\to 0}  \Big( C(r)+F(r)+ D_*(r)\Big) \ < \ \ep_*^{\frac 13}.
$$
Theorem \ref{CKN_condition} is proved.

\bigskip
\noindent
{\bf Proof of Theorem \ref{CKN_theorem}:}
Assume condition (\ref{ep-regularity-1}) holds and let $\ep_0(K)>0$ be the constant defined in Theorem \ref{CKN_condition}.
From (\ref{ep-regularity-1}), (\ref{ep-regularity-2}) we obtain there exists $R>0$ such that
$$
\sup\limits_{r\in (0,R)} E(r) \ < \ \ep_0 \qquad\mbox{and}\qquad  \sup\limits_{r\in (0,R)} E_*(r) \ < \ K.
$$
Denote $(v^R,H^R,p^R)$ the functions
$$
\gathered
v^R(x,t)=R v(x_0+ Rx,t_0+R^2t), \qquad H^R(x,t)=R H(x_0+Rx,t_0+R^2t), \\ p^R(x,t)=R^2 p(x_0+Rx,t_0+R^2t).
\endgathered
$$
Then functions $(v^R,H^R,p^R)$ satisfy all conditions of Theorem \ref{CKN_condition}. The result follows from Theorem \ref{CKN_condition} and Theorem \ref{Fixed_r}.

\bigskip
\noindent
{\bf Proof of Theorem \ref{Partial_Regularity}:} The result is a direct consequence of  Theorem \ref{CKN_theorem} and measure theory, see \cite{CKN}, \cite{Lin}, \cite{LS}, \cite{Seregin_JMFM}.

\end{document}